%% file: RealOption.main.tex
\documentclass[11pt]{amsart}
\usepackage[T1]{fontenc}
\usepackage[utf8]{inputenc}
\usepackage[english]{babel}
\usepackage{amsmath, amsthm, amssymb, amsfonts, enumerate, mathrsfs, dsfont, comment, bm, mathtools}
\usepackage{algorithm}
\usepackage{algpseudocode}
\usepackage{enumitem}
\usepackage[colorlinks=true,linkcolor=blue,urlcolor=blue]{hyperref}
\usepackage{color, geometry, todonotes, indentfirst, graphicx, pdfpages}
\usepackage[title]{appendix}
\usepackage{fancyhdr}
\usepackage{polski}
\allowdisplaybreaks

\newtheorem{theorem}{Theorem}[section]
\newtheorem{remark}[theorem]{Remark}
\newtheorem{assumption}[theorem]{Assumption}

\newtheorem{proposition}[theorem]{Proposition}

\theoremstyle{plain}
\def \dd{{\rm d}}
\def \R{\mathbb{R}}

\def \F{\mathbb{F}}

\def \cF{{\mathcal F}}

\def \A{\mathcal{A}}

\def \and{\quad \text{and} \quad}

\title[RL for Exploratory Real Option]{Reinforcement Learning in Real Option Models}
\author[Dianetti]{Jodi Dianetti}
\author[Ferrari]{Giorgio Ferrari}
\author[Xu]{Renyuan Xu}

\date{ \today. \\{\it Contacts:} \href{jodi.dianetti@uniroma2.it }{jodi.dianetti@uniroma2.it}, \href{giorgio.ferrari@uni-bielefeld.de }{giorgio.ferrari@uni-bielefeld.de},  \href{renyuanxu@stanford.edu}{renyuanxu@stanford.edu}}

\numberwithin{equation}{section}

\begin{document}

\begin{abstract}
We investigate an entropy-regularized reinforcement learning (RL) approach to optimal stopping problems motivated by real option models. Classical stopping rules are strict and non-randomized, limiting natural exploration in RL settings. To address this, we introduce entropy regularization, allowing randomized stopping policies that balance exploitation and exploration. We derive an explicit analytical solution to the regularized problem and prove convergence of the associated free boundary to the classical stopping threshold as the entropy vanishes. The regularized problem admits a natural formulation as a singular stochastic control problem. Building on this structure, we propose both model-based and model-free policy iteration algorithms to learn the optimal boundary. The model-free method operates without knowledge of system dynamics, using only trajectories from the stochastic environment. We establish convergence guarantees and illustrate strong numerical performance. This framework provides a principled and tractable approach for data-driven stopping problems under uncertainty.
\end{abstract} 
\maketitle  




\section{Introduction} 
\label{sec:intro}

Optimal stopping problems are foundational in the fields of decision theory, economics, finance, and operations research. At their core, these problems concern determining the most advantageous time to take a particular action in order to maximize expected rewards or minimize costs under uncertainty. Among the most widely known applications are real options, where the timing of investments, abandonments, or exits in uncertain environments plays a critical role. Classical approaches to such problems rely heavily on precise modeling of stochastic processes and the derivation of analytic or semi-analytic solutions through variational inequalities and free boundary characterizations (we refer to \cite{DixitPindyck, peskir.shiryav.2006optimal, stokey2008economics} for methodologies and examples in Economics and Finance).

However, in real-world applications, decision-makers often face limited knowledge about the underlying dynamics governing the system or incomplete access to critical model parameters. This lack of model fidelity makes the practical implementation of optimal stopping policies significantly more complex. Recent advances in data-driven decision-making, particularly through reinforcement learning (RL), offer an avenue to overcome such model uncertainty by learning optimal or near-optimal policies directly from simulated or observed trajectories, circumventing the need for full specification of the system's dynamics.

This paper contributes to the growing intersection of optimal stopping theory and reinforcement learning by studying a continuous-time real option problem through the lens of entropy-regularized singular stochastic control. Real option problems arise in a wide range of economically significant contexts, including irreversible investment \cite{bulan2005real,bulan2005real,guo2005irreversible}, project abandonment \cite{cruz2016assessing,vintila2007real}, capacity expansion \cite{bensoussan2019sequential}, technology adoption \cite{huisman2004strategic,alvarez2001adoption,xu2023decision}, and resource extraction \cite{mezey2010real,trigeorgis1996real}. Understanding how and when to exercise such options under uncertainty is a central question in corporate finance, macroeconomics, and operations management.

Our approach departs from classical formulations by embedding the stopping decision within a family of singular controls and augmenting the reward functional with an entropy regularization term. This regularization, inspired by the reinforcement learning literature, provides a systematic way to incorporate exploration into the decision process, thereby offering a new perspective on the structure and robustness of optimal exercise policies. While this paper focuses on understanding the advantages of entropy regularization and on exploiting the singular control formulation by constructing an explicitly solvable solution, a general theoretical analysis of entropy-regularized optimal stopping problems is developed in our companion paper \cite{dianetti.ferrari.xu.2024exploratory}. Importantly, the assumptions in that work do not cover the class of real option problems considered here.


The classical real option problem considered here involves a stochastic process $X^x$ governed by a geometric Brownian motion, where a decision-maker accrues running profits modeled by a concave function $\pi(\cdot)$ and faces a discounted lump-sum reward $\kappa$ upon exit. The objective is to determine the optimal stopping time $\tau^*$ that maximizes the expected discounted cumulative payoff. Under regularity and monotonicity assumptions, the value function and optimal policy are well-understood: optimal stopping is triggered when the state variable $X^x$ hits a problem-dependent threshold, and the value function satisfies a free-boundary problem that can be solved analytically or semi-analytically.

Despite the tractability of this classical solution, it inherently presumes full knowledge of the model parameters (e.g., drift $\mu$, volatility $\sigma$, discount rate $\rho$, and profit function $\pi$). In practice, these parameters are often estimated imprecisely or evolve over time, rendering such solutions suboptimal or unusable. Moreover, the deterministic nature of the optimal stopping rule (a sharp threshold) leads to non-exploratory behavior, which limits the agent's ability to learn the optimal policy when operating under uncertainty.

One key challenge that distinguishes OS from regular control is that optimal stopping involves a non-smooth decision -- whether to stop or continue -- while regular controls gradually change the dynamics through drift and/or volatility coefficients. 
As a result, gradient-based RL algorithms, although popular for regular controls, cannot be directly applied to the ``stop-or-continue'' decisions in our setting \cite{reppen2022neural,soner2023stopping}. 
To overcome this difficulty, we replace these sharp rules based on the hitting times by stopping probabilities, leading to a {\it randomized stopping time} -- encoded as singular controls -- that follows the stopping probability. 
Unlike the heuristic of a fuzzy boundary considered in \cite{reppen2022neural}, where the agent stops with a probability proportional to the distance to the boundary, our approach to randomized stopping time is based on principled guidance -- regularized objective function with cumulative residual entropy.  
Another advantage of randomized stopping time is its exploratory nature; it stops at different scenarios according to certain probability, thereby collecting more information from the unknown environment in the RL regime. 
Compared to classic RL settings, exploration is particularly essential for optimal stopping problems because the terminal reward can only be collected upon making the stopping decision, contributing to the {\it reward sparsity} challenge in RL \cite{devidze2022exploration,hare2019dealing}.
Therefore, randomized stopping time helps the decision-maker learn more about the unknown terminal reward.

Our contributions are both theoretical and algorithmic. First, we derive a closed-form solution to the entropy-regularized control problem. 
This takes the form of a two-dimensional degenerate singular stochastic control problem (see, e.g., \cite{angelis&ferrari&moriarty2019, ferrari2015integral, guo&tomecek2009}), for which we are able to provide the explicit solution by solving the Hamilton-Jacobi-Bellman (HJB) variational inequality with gradient constraint associated to the value function. The optimal randomized stopping time takes the form of a reflection policy at a state-dependent reflection boundary. This is a continuously differentiable function of the "fuel" variable, representing the probability of not having stopped up to a certain time. The explicit nature of our results allows also to rigorously establish the convergence of the resulting free-boundary to the classical stopping boundary as the entropy weight vanishes. 

Second, we develop a reinforcement learning algorithm -- based on policy iteration -- that learns the optimal boundary in both model-based and model-free settings. The algorithm iteratively updates the reflection threshold by solving a parameter-dependent ODE in the exploration region and adjusting the policy based on gradient information. 
 Importantly, we provide theoretical guarantees on policy improvement and, under suitable initialization, on convergence of the learned boundary to the entropy-regularized optimum. {Policy improvement results are more commonly established in the continuous-time RL literature, whereas proving policy convergence remains highly nontrivial, with only a few notable exceptions in standard control settings \cite{bai2023reinforcement, huang2022convergence, ma2024convergence, tran2024policy}. To date, no such convergence results have been established for optimal stopping or singular control problems prior to our work. A key challenge lies in establishing the regularity of the free boundary throughout the entire iterative procedure. Our approach draws inspiration from \cite{kumar2004numerical} on a one-dimensional singular control problem, but their convergence proof hinges on a second-order condition that must hold at every iteration—a requirement that is both strong and hard to verify.} In the model-free version, we demonstrate how the boundary can be learned through trajectory-based sampling and simulator interaction, without knowledge of model parameters or reward function structure.

To the best of our knowledge, this is the first work to combine singular control, entropy regularization, and reinforcement learning in a continuous-time real option setting. Our theoretical analysis bridges optimal stopping, singular stochastic control, and modern RL. From an algorithmic perspective, our policy iteration scheme for boundary learning sets the stage for scalable extensions to high-dimensional and path-dependent settings.\\

\paragraph{\bf Related literature.} Reinforcement learning (RL) for optimal stopping problems is intimately connected to the broader challenge of RL with sparse rewards \cite{devidze2022exploration,hare2019dealing}. In this context, the reward $G(X_\tau)$ is received only at the stopping time, which leads to extreme reward sparsity. This sparsity introduces substantial learning difficulties when compared to more standard control tasks -- whether in continuous time or classical discrete-time Markov decision processes.

In settings where the underlying model is fully specified, \cite{herrera2021optimal} proposed using randomized neural networks to approximate solutions of optimal stopping problems, with applicability to high-dimensional settings.  Recent works such as \cite{reppen2022neural} and \cite{soner2023stopping} have proposed deep learning algorithms to approximate optimal stopping boundaries. In a related direction, \cite{bayer2021randomized} analyzed randomized optimal stopping problems and established convergence results for both forward and backward Monte Carlo-based optimization methods. Additionally, \cite{ata2023singular} demonstrated the efficacy of deep learning approaches in solving high-dimensional singular control problems.

A notable contribution to continuous-time RL comes from \cite{denkert2024control}, who developed a general framework for policy gradient methods that applies to a wide class of control problems, including optimal stopping, impulse, and switching problems. Their method exploits the connection between stochastic control and randomized formulations. However, this framework currently lacks theoretical convergence guarantees.

The closest work to ours is \cite{dong2024randomized}, which investigates a regularized version of the American Put Option under a Shannon entropy framework in one dimension. The study adopts an intensity control formulation to induce exploration and proves the convergence of the policy iteration algorithm (PIA) under fixed temperature parameters. Nevertheless, the convergence of the optimal policy as the temperature $\lambda \to 0$ remains unresolved—even in this simple setting.

Our work departs from \cite{dong2024randomized} in several key ways. First, we employ a regularity-based analysis that combines partial differential equation (PDE) techniques with a probabilistic link to another optimal stopping problem. In contrast, the analysis in \cite{dong2024randomized} relies heavily on explicit computations, which are less adaptable to high-dimensional settings. Second, we establish the convergence of the optimal policy for the regularized problem to the original optimal stopping policy as $\lambda \to 0$, irrespective of the problem's dimension—something not achieved in \cite{dong2024randomized}.

Recently, \cite{dai2024learning} proposed an alternative method by reformulating the optimal stopping problem as a penalized regular control task with binary decisions (0 or 1). Entropy regularization is then applied, converting the setup into a classic entropy-regularized stochastic control problem, to which standard RL algorithms can be applied.

Another important line of research uses non-parametric statistics for stochastic processes to develop learning-based control algorithms \cite{christensen2023data2,christensen2023nonparametric,christensen2024learning,christensen2023data}. Specifically, \cite{christensen2023nonparametric} and \cite{christensen2024learning} tackled one-dimensional singular and impulse control problems by learning critical value thresholds within a non-parametric diffusion framework. This methodology was extended to multivariate reflection problems in \cite{christensen2023data}, while \cite{christensen2023data2} focused on theoretical guarantees such as regret bounds and non-asymptotic PAC estimates.

The rest of the paper is organized as follows. In Section~\ref{sec:real_option}, we present the classical and entropy-regularized formulations of the real option problem, including the singular control reformulation and the entropy-based exploration term. Section~\ref{sec:real_option_analytical}  derives the explicit analytical solution to the entropy-regularized problem and studies its limiting behavior. Section~\ref{section 5} develops the reinforcement learning algorithms for both model-based and model-free scenarios, including implementation details, convergence results, and numerical experiments. Finally, concluding remarks are made in Section \ref{sec:conclusion}, while proofs are collected in Section \ref{Appendix A}.


\section{The entropy-regularized real option problem}
\label{sec:real_option}

Let $(\Omega, \mathcal{F}, \mathbb{P})$ be a complete probability space. We assume that the probability space is large enough to accommodate a one-dimensional Brownian motion $W:=(W_t)_{t\geq0}$ and a $W$-independent uniformly distributed random variable $U:\Omega \to [0,1]$. We denote by $\mathbb{F}^W:=(\mathcal{F}^W_t)_{t\geq0}$ the filtration generated by $W$, as usual augmented by the $\mathbb{P}$-null sets of $\mathcal{F}$. 

For a given running profit $\pi:\mathbb{R}_+\rightarrow\mathbb{R}$, consider the classical real option/optimal stopping (OS) problem \cite{DixitPindyck}
\begin{equation}
\label{eq:OS0}
\begin{aligned}
& \sup _{\tau \in \mathcal T} \mathbb{E}\left[\int_{0}^{\tau} e^{- \rho s}  \pi(X^x_{s}) \dd  s+\kappa e^{-\rho \tau}\right], \\
 & \text{subject to }  \dd X^x_{t}  =\mu X^x_t \dd t +\sigma X^x_t \dd W_{t}, \quad X^x_{0}=x>0.
\end{aligned}
\end{equation}
Here, a decision  maker aims at determining the best time to leave the market in exchange of the lump-sum payoff $\kappa>0$. Up to that exit time, stochastic profits are accumulated and discounted at rate $\rho>0$. Uncertainty in the instantaneous returns of the company is given through the process $X^x$, which models a reference economic indicator, grows at exponential rate $\mu \in \R$, and it is subject to Brownian fluctuations with volatility $\sigma>0$. 

Noticing that, 
$$
 \mathbb{E}\left[\int_{0}^{\tau} e^{- \rho s}  \pi(X^x_{s}) \dd s+\kappa e^{-\rho \tau}\right]
 = \mathbb{E} \left[\int_{0}^{\tau} e^{-\rho s}\left( \pi(X^x_{s})- \rho \kappa \right) \dd s\right]+\kappa, 
$$
we therefore consider from now on the OS problem
\begin{equation}\label{eq real option OS}
V(x) : =\sup _{\tau \in \mathcal T} \mathbb{E} \left[\int_{0}^{\tau} e^{-\rho s}\left( \pi(X^x_{s})- \rho \kappa\right) \dd s\right].
\end{equation}

We make the following standing assumption.
\begin{assumption}\label{ass:Pi}
     The data of the problem satisfy the following conditions:
    \begin{enumerate}
        \item $\pi\in{C}^2(\mathbb{R}_+)$ with  $\pi(0)=0$, $\pi$ is non-decreasing and concave, $x\pi'(x)$ is strictly increasing and $\pi$ has a sublinear growth; i.e., there exists some $\theta\in(0,1)$ and $c>0$ such that $$|\pi(x)|\leq c(1+|x|^\theta).$$
        \item $\rho$ is large enough: namely,  $\rho >\mu$. 
    \end{enumerate}
\end{assumption}
\noindent One relevant example that satisfies Assumption \ref{ass:Pi} is $\pi(x)= x^{\theta}$ for $\theta \in (0,1)$. Notice that, under Assumption \ref{ass:Pi}, the value in \eqref{eq real option OS} is finite.

When the data of the model $\mu, \sigma,\pi,\rho,\kappa$ are known,
Problem \eqref{eq real option OS} can be explicitly solved via classical arguments (see, e.g., \cite{DixitPindyck}) based on the application of the smooth-fit and smooth-pasting principles and a verification argument. In particular, it turns out that the optimal stopping time is given by
$$\tau^*:=\inf\{t\geq0:\, X^{x}_t \leq b^*\},$$
for some free boundary $b^*>0$ that is the unique solution to a certain algebraic equation depending on the problem's data.

This is clearly not our interest. We are concerned here in illustrating the relations between such a classical solution approach and the entropy regularization approach (see in particular the result of Theorem \ref{theorem RO free boundary convergence} below) and in developing RL algorithms able to learn the optimal stopping rule, specifically in model-free settings in which the decision maker does not have access to precise evaluations of the problem's data. This will be achieved by introducing an entropy-regularized version of the previous real option problem, as discussed in the folllowing sections.

\subsection{Exploratory formulation via singular controls}

We briefly recap here the approach introduced in the companion paper \cite{dianetti.ferrari.xu.2024exploratory}, where the theoretical framework for RL algorithms for multi-dimensional OS problems is provided.

Inspired by \cite{touzi.vieille.2002continuous}, we introduce a notion of randomized/exploratory stopping times. For any $y\in[0,1]$, define the set of stochastic processes
\begin{align}
\label{eq:setAy}
    \A (y):= \Big\{ \xi : & \, \Omega \times[0, \infty) \to [0,y], \text{ $\F ^W$-adapted, nondecreasing, c\`adl\`ag, with $\xi_{0-}=0$} \Big\}. 
\end{align} 
In the rest of this paper, an element of $\A (y)$, $y\in [0,1]$, will be referred to as \emph{singular control}.

Given any $\xi \in \mathcal A (1)$, we consider the random variable $\tau^\xi$ on $(\Omega, \mathcal F)$ by setting $\tau^\xi:=\inf \left\{t \geq 0 \, | \, \xi_{t}>U\right\}$, with the convention $\inf \emptyset = +\infty$, and we will refer to it as \emph{randomized/exploratory stopping time}.
Notice that $\tau^\xi$ is not necessarily an $\mathbb F ^W$-stopping time.
However, if $\tau \in \mathcal T$, the process $\xi^\tau \in \mathcal A (1)$ defined by $\xi^\tau := ( \mathds 1 _{\{ t \geq \tau \} } )_t$ is such that $\tau = \tau^{\xi^\tau}$, which gives a natural inclusion of $\mathcal T$ into $\mathcal A (1)$.
Moreover, we have
$$
\mathbb{P} \big( \tau^\xi \leq t \mid \cF _{t}^{W}\big)
=\mathbb{P}\left( U \leq {\xi}_{t} \mid \cF _{t}^{W}\right)
=\int_{0}^{\xi_{t}} \dd u=\xi_{t}, 
$$
so that $\xi_t$ can be interpreted as the probability of stopping before time $t$.
Notice also that $\xi_\infty := \lim_{t \to \infty} \xi_t$ is not necessarily equal to 1: indeed, the related randomized stopping time $\tau^\xi$ is not necessarily finite.

For $\xi \in \mathcal A (1)$, by evaluating the performance of the randomized stopping time $\tau^\xi$, we obtain
\begin{align*}
J(x;\tau^\xi)  =& 
\mathbb{E}\left[ \int_{0}^{\tau^\xi} e^{- \rho t} \big(\pi (X^x_{t}) - \rho \kappa\big) \dd t \right] 
 = \mathbb{E}\left[\int_{0}^{\infty} e^{- \rho t} \big(\pi (X^x_{t}) - \rho \kappa\big) \mathds{1}_{\{t \leq \tau^\xi \}} \dd t \right] \\
 =& \mathbb{E}\left[\int_{0}^{\infty} e^{-\rho t} \big(\pi (X^x_{t}) - \rho \kappa\big)  \mathds{1}_{\left\{\xi_{t}< U\right\}} \dd t \right] 
 = \mathbb{E}\left[\int_{0}^{\infty} e^{-\rho t} \big(\pi (X^x_{t}) - \rho \kappa\big)\left(1-\xi_{t}\right) \dd t \right].
\end{align*}
This suggests to define a profit functional in terms of the singular controls $\xi \in \mathcal A (1)$ by
\begin{equation*}
    J^0 (x;\xi)  := \mathbb{E}\left[\int_{0}^{\infty} e^{-\rho t} \big(\pi (X^x_{t}) - \rho \kappa\big) \left(1-\xi_{t}\right) \dd t \right], \quad \xi \in \mathcal A (1),
\end{equation*}
as a natural extension of $J$ to randomized stopping times.

However the new payoff functional is linear in the singular control $\xi$. As a consequence, as expected, allowing for randomized stopping times does not change the optimal value and (crucial for our story) it  does not imply that optimal actions are necessarily randomized. This is formally proved in the next proposition (whose proof is postponed to Section \ref{sec:AppA}) and the subsequent remark.

\begin{proposition}\label{prop xi does not explore}
For any $x \in \R^n$ we have 
$$
V(x) = \sup _{\xi \in \mathcal A (1) } J^0(x;\xi) = J^0 (x;\xi^*), 
$$
with $\xi^* := ( \mathds 1 _{ \{ t \geq \tau^* \} } )_t$ and $\tau^*$ optimal for $J(x;\cdot)$.    
Moreover, $\xi^* := (\mathds 1 _{ \{ t \geq \tau^* \} })_t$ is the unique optimal control for $J^0(x;\cdot).$
\end{proposition}

\begin{remark}[Non-exploratory behavior of the optimal controls]\label{remark xi does not explore}
We revisit Proposition \ref{prop xi does not explore} from both reinforcement learning (RL) and exploratory perspectives. The singular control \( \xi^* \) associated to the optimal stopping time $\tau^*$ aligns precisely with this strict stopping time. Consequently, the optimal strategy does not require any form of randomization. In effect, the agent halts the process as soon as an action becomes necessary, which precludes any gradual gathering of information about the performance of alternative actions.

From an RL standpoint, this non-exploratory nature of the optimal control \( \xi^* \) highlights a disconnect between optimization and information acquisition. It implies that the agent, in pursuing the optimal strategy, does not naturally explore other possibilities. Therefore, some modification is required to address this limitation. This is addressed in the next section.
\end{remark}

\subsection{Entropy regularization}
In light of Proposition \ref{prop xi does not explore} and of the related Remark \ref{remark xi does not explore}, we introduce a regularization term to incentivize exploration/randomization.  This is achieved by incorporating an entropy term to regularize the problem, drawing motivation from the RL literature.

Since the sample paths of the processes $\xi \in \mathcal A (1)$ are not necessarily absolutely continuous with respect to the Lebesgue measure, we choose (a discounted version of) the \emph{cumulative residual entropy} (see \cite{rao.chen.vemuri.wang.2004cumulative}) weighted by a parameter $\lambda >0$; 
namely, we consider, for $\xi \in \A (1)$ the entropy
\begin{equation}\label{eq:CRE}
\begin{aligned}
\Lambda^\lambda (\tau^\xi) :=&-\lambda \int_{0}^{\infty} e^{- \rho t} \mathbb{P} \left( \tau^\xi \geq t \mid \cF _{t}^{W} \right) \log \left(\mathbb{P}\left(\tau^\xi \geq t \mid \mathcal{F}_{t}^{W}\right)\right)  \dd t \\
 =& -\lambda \int_{0}^{\infty} e^{- \rho t}\left(1-\xi_{t}\right) \log \left(1-\xi_{t}\right)  \dd t =: \Lambda^\lambda (\xi).
\end{aligned}
\end{equation}
The entropy $\Lambda^\lambda$ is nonnegative and achieves its highest values when the probability $\xi_t$ is near the level $e^{-1}$, thus incentivizing the use of randomized stopping times.

Building on this intuition,  we define an entropy regularized OS problem 
\begin{equation}\label{eq value singular control problem}
    V^\lambda (x) := \sup _{\xi \in \mathcal A (1)} J^\lambda ( x; \xi ), 
\end{equation}
where the exploration--exploitation trade off is captured by the profit functional 
\begin{equation*}
    J^\lambda ( x; \xi ):= \mathbb E \bigg[ \int _0^\infty e ^{-\rho t} \underbrace{\big( \pi(X^x_t) - \rho \kappa\big) (1-\xi_t) }_{\text{exploitation}}\dd t  -\lambda \int_{0}^{\infty} e^{- \rho t} \underbrace {\left(1-\xi_{t}\right) \log \left( 1-\xi_{t} \right)}_{\text{exploration}} \dd t \bigg].
\end{equation*}
A control $\xi^\lambda \in \mathcal A (1)$ is said to be optimal if $J^\lambda (x;\xi^\lambda ) = V^\lambda (x)$.

In order to embed the new optimization problem into a Markovian framework, we introduce the extra state variable $y \in [0,1]$ and, for $\lambda >0$, define the singular control problem
\begin{equation}\label{eq entropy real option OS}
\begin{aligned}
V^\lambda (x, y) :=&\sup _{\xi \in \A(y)} \mathbb{E}  \left[ \int_0^{\infty}  e^{- \rho t} \big( \pi  (X^x_{t} )- \rho \kappa ) Y_{t}^{y,\xi}-\lambda Y_{t}^{y,\xi} \log  (Y_{t}^{y,\xi} ) \big) \dd t \right], \\
 \text{subject to }\,\,  \dd X^x_{t} =&\mu X^x_t \dd t +\sigma X^x_t \dd W_{t}, \quad  X^x_{0}=x, \\
 \dd Y^{y,\xi}_t =& -\dd\xi_t, \quad Y^{y,\xi}_{0-}=y, 
\end{aligned}
\end{equation}
with $\A(y)$ as defined in \eqref{eq:setAy}.

Clearly (with a slight abuse of notation), we have by \eqref{eq value singular control problem} and \eqref{eq entropy real option OS} that $V^{\lambda}(x)= V^{\lambda}(x, 1)$.

Standard arguments based on the dynamic programming principle leads to the HJB equation expected to be solved (in suitable sense) by $V^\lambda$:
\begin{equation*}
    \max \Big\{ (\mathcal{L}_x-\rho) V^\lambda + ( \pi(x)-\rho \kappa ) y -\lambda y \log y, -V^\lambda_{y} \Big\}=0, 
    \end{equation*}
with boundary condition $V(x, 0)=0$, $x \in (0,\infty)$, where  $\mathcal L _x \phi (x) := \mu x \phi_x (x) + \frac{1}{2} \sigma^2 {x^2} \phi_{xx} (x) $.


\section{The analytical solution to the entropy-regularized real option problem}
\label{sec:real_option_analytical} 

In this section, we provide the explicit solutions to the entropy-regularized real option problem \eqref{eq entropy real option OS}. The proof of the next theorem, which is the main result of this section, is postponed to Section \ref{sec:AppB}.

\begin{theorem}[The solution to the entropy-regularized real option problem]
\label{theorem verification real option}
Let $(x,y) \in (0,\infty) \times [0,1]$ be given and fixed.
Introduce the nondecreasing function
\begin{eqnarray}\label{eq:g_lam}
    g_{\lambda}(x) := \exp \left( \frac{-\frac{H_\pi'(x)}{\alpha_-}x+H_\pi(x)-\kappa - \frac{\lambda}{\rho}}{\frac{\lambda}{\rho}}\right) \wedge 1,
\end{eqnarray}
such that $y^\lambda:= g_\lambda(0) = e^{-1-\frac{\kappa \rho}{\lambda} }$,
and define
$$H_\pi(x) := \mathbb{E}\left[\int_0^\infty e^{-\rho t}\pi(X_t^x)\dd t\right], \quad A_{2}(y) :=\int_{g_\lambda(0)}^y  \frac{\kappa  +\frac{\lambda}{\rho}  \log (u)+ \frac{\lambda}{\rho}- H_\pi(g_\lambda^{-1}(u))}{(g_\lambda^{-1}(u))^{\alpha_-}} \dd u.$$

Then, letting $F(x, y):=A_{2} (y) \,x^{\alpha_-} + H_\pi(x)y - \kappa y -\frac{\lambda}{\rho} y \log y$, the value of the entropy-regularized real option problem \eqref{eq entropy real option OS} is given by 
\begin{equation}
\label{valuefunction-thm}
V^{\lambda}(x,y) =
\begin{cases}
H_\pi(x)y - \kappa y - \frac{\lambda}{\rho} y \log y, & y<y^\lambda, \\
F(x, y), & y\leq g_\lambda (x), \ y \geq y^\lambda, \\
F(x, b_\lambda^{-1} (x) ), & y> g_\lambda (x), \ y \geq y^\lambda, 
\end{cases}
\end{equation}
and the reflection policy 
$$\xi_{t}^\lambda := 
\sup\limits_{0 \leq s \leq t} \big( y - g_{\lambda} (X_s^x) \big)^{+}, \quad t>0, \quad \xi_{0^-}^\lambda =0,$$ is optimal for the initial condition $(x,y) \in (0,\infty) \times [0,1]$.
\end{theorem}


Given the explicit quantification of the solution in Theorem \ref{theorem verification real option}, we can now determine  the vanishing entropy limit of the optimal reflection boundary. This will agree with the optimal stopping boundary of the classical real option problem without entropy-regularization.

Recall that the standard OS problem \eqref{eq real option OS} can be solved semi-explicitly (see, e.g., Chapter 5 in \cite{pham2009continuous} or \cite{DixitPindyck}) via a classical guess-and-verify approach and an application of the smooth-fit principle. 
In particular, the optimal stopping time for such a problem is given by
$$
\tau^{*}=\inf \left\{t \geq 0: X_{t}^{x} \leq b^{*}\right\}, 
$$
with $b^*$ defined as the solution to the free boundary equation: 
\begin{equation}\label{eq RO free boundary equation reflecting}
    \frac{H'_\pi(x)}{- \alpha_-}x +H_\pi(x)= \kappa.
\end{equation}
Indeed, the monotonicity conditions on $\pi$ in Assumption \ref{ass:Pi} ensure that the solution to \eqref{eq RO free boundary equation reflecting} is unique.
For example, in the special case of $\pi(x) = x^\theta$ with $\theta\in(0,1)$, we would have $b^*=\left[-\frac{1}{P}\left(\frac{\alpha_-}{\theta-\alpha_-} \right) \kappa \right]^{1 / \theta}$ with $P := \frac{1}{\rho+\frac{1}{2} \sigma^{2} \theta(1-\theta)-\theta \mu}$.
\\ \indent
We next look at the entropy regularized OS problem.  
Notice that, for $y^\lambda$ as in \eqref{eq definition of minimal y}, we  have the limits  $ \lim _{\lambda \downarrow 0} y^\lambda = 0$.
Moreover, the free boundary of the entropy regularized OS problem can be equivalently expressed in terms of the function  $b_\lambda:=g_\lambda^{-1}: [y^\lambda, 1] \to [0,\infty)$. 
By \eqref{eq:optimal_g}, we deduce that $b_\lambda$ must satisfy the equation
\begin{equation}\label{eq free boundary limits}
  \frac{H_\pi'(b_\lambda(y))}{ - \alpha_-}b_\lambda(y)\,+\,H_\pi(b_\lambda(y))= \kappa + \frac{\lambda}{\rho} (1+\log(y)). 
\end{equation}
For $y\in[e^{-1},1]$, we have $1+\log(y) \leq 0$. 
Thus, since   $(H'_\pi(x)x)'>0$ and $H_\pi'\ge 0$ (according to Assumption \ref{ass:Pi}), it follows that $b_\lambda (y) \geq b^*$ and that the map $\lambda \mapsto b_\lambda(y)$ is nonincreasing.
Similarly, for $y \in (0,e^{-1}]$,  we have $1+\log(y) \geq 0$, so that $b_\lambda (y) \leq b^*$ and  the map $\lambda \mapsto b_\lambda(y)$ is nondecreasing.
Hence, for any $y \in (0,1]$, there exists a finite limit $b_0 (y) := \lim_{\lambda \to 0} b _\lambda (y)$. 
Taking limits in \eqref{eq free boundary limits} as $\lambda \to 0$, we obtain that, for any $y \in (0,1]$,  $b_0(y)$ is a solution to \eqref{eq RO free boundary equation reflecting}. 
By uniqueness of the solution to \eqref{eq RO free boundary equation reflecting}, we obtain that $b_0(y) = b^*$ for any $y \in (0,1]$.
\\ \indent
Thus, we have proved the following result.
\begin{theorem}[Vanishing entropy limits]
    \label{theorem RO free boundary convergence}
    The optimal reflection boundary $b_\lambda$ of the entropy regularized OS problem \eqref{eq entropy real option OS} converges to the  optimal stopping boundary $b^*$ of the OS problem \eqref{eq real option OS}; namely, 
    $\lim_{\lambda \to 0} b_\lambda (y) = b^*$, for any $y \in (0,1]$.
\end{theorem}


\section{The reinforcement learning algorithm}
\label{section 5}

In this section, we propose an RL framework to learn the reflection boundary  $g_\lambda$ and the corresponding value function $V^\lambda$, without prior knowledge on the form of $g_\lambda$. Our proposed RL framework has  two versions: a model-based numerical version (see Subsection \ref{sec:model-based}) and a model-free learning version (see Subsection \ref{sec:model-free}). In the first version, where all model parameters are known, we design a Policy Iteration algorithm to numerically find $g_\lambda$. In the second version, where the model parameters are unknown, we provide a sample-based Policy Iteration method to learn $g_\lambda$.

While our RL framework will be formulated in terms of the real option example, we believe the results can be extended to a more general setting and in higher dimensions.

 In line with RL literature, from now on the value function $V^\lambda$ will be referred to as \textit{optimal value function}, while the profit $J^\lambda (x,y;\xi)$ of a given policy $\xi \in \mathcal A (y)$ will be referred to as  \textit{value function} associated to the policy $\xi$.
Moreover, we set $\R_+ := [0,\infty)$.

\subsection{Model-based numerical analysis}\label{sec:model-based}
In this subsection, we present a numerical scheme to solve for the optimal boundary, utilizing full knowledge of the underlying system. As such, this method is termed a model-based analysis.

Motivated by the previous results (see Theorem \ref{theorem verification real option}, here we focus on a subclass of control policies that can be fully characterized by some reflection boundary $g$ that satisfies the following assumption.
\begin{assumption}\label{ass:g}
Assume $g:\mathbb{R}_{+}\rightarrow (0,1]$ is a non-decreasing function such that     $$g\in {C}^1([0,\hat{x}_g]), \quad \text{where} \quad \hat{x}_g = \inf_{x\in \mathbb{R}_+}\{g(x)=1\}.$$
  \end{assumption}
For a function $g$ that satisfies Assumption \ref{ass:g}, we define the associated policy $\xi^g$  of the following form:
  \begin{eqnarray}\label{eq:parameterized_policy}
      \xi^g_t = \sup_{s\leq t} \big(y-g(X_s^x)\big)^+, t \ge 0.
  \end{eqnarray}

Define the associated value function as: 
 \begin{equation} 
 \begin{aligned}
      V^{\lambda}_g (x,y) := &J^{\lambda}(x,y;\xi^g) \nonumber \\
      =&\mathbb{E}\left[ \int_0^\infty e^{-\rho t} \Big(\Big( \pi(X^x_t)-\rho \kappa\Big)Y^{y,\xi^g}_t - \lambda Y^{y,\xi^g}_t \log (Y^{y,\xi^g}_t)\Big) \dd t \right]
  \end{aligned}
  \end{equation}
  subject to $Y^{y,\xi^g}_t = y-\xi_t^g$, with $\xi_t^g$ defined in \eqref{eq:parameterized_policy}.

The policy in \eqref{eq:parameterized_policy} defines two areas: the exploration area and the stopping area, both associated with $g$:
  \begin{eqnarray}
      \mathcal{E}(g) &:=&\Big\{(x,y)\,\,\Big|\,\, y \leq g(x) \Big\}, \\
      \mathcal{S}(g) &:=& \Big\{(x,y)\,\,\Big|\,\,   y> g(x)  \Big\}.
  \end{eqnarray}

Then we have the following results for the value function $V^\lambda_g$ associated with policy $\xi^g$. Its proof is deferred to  Section \ref{sec:AppC}.
 \begin{theorem}\label{thm:property_g}
    Assume Assumptions \ref{ass:Pi} and \ref{ass:g} hold. Then the following  boundary value problem has a unique ${C}^{1}(\mathbb{R}_+\times [0,1])\cap{C}^{2}\big(\overline{\mathcal{E}(g)}\big)$ solution:
     \begin{eqnarray}
      &&(\mathcal{L}_x-\rho) u + \big(\pi(x)-\rho \kappa \big)y - \lambda y \log y                   =0\quad \text{on} \quad\mathcal{E}(g),\\
     && -u_y = 0, \quad \text{on}\quad \mathcal{S}(g).
  \end{eqnarray}
  
In addition, the solution satisfies $u_{xx}\in \mathbb{L}_{loc}^{\infty}(\mathbb{R}_+\times[0,1])$. 
Finally, we have   
\begin{eqnarray}
   V^\lambda_g(x,y) \equiv u(x,y). \label{eq:verification_thm}
\end{eqnarray}
\end{theorem}

With the results in Theorem \ref{thm:property_g}, we update the new boundary following:
  \begin{equation}\label{eq:updated_policy}
    \widetilde{g}(x) =
      \begin{cases}
       & \max \Big\{y< g(x)\,\Big| \partial_{xy} V^\lambda_{g}(x,y) = 0 \Big\} \,\,\text{if} \,\, \partial_{xy}^- V^\lambda_{g}(x,g(x))<0, \text{ and } \\
    &   \widetilde{g}(x) = g (x)  \qquad \text{otherwise,}  
      \end{cases}
  \end{equation}
where we define the notation $\partial_{xy}^- f(x,y):= \lim_{h\rightarrow0-}\frac{\partial_x f(x,y+h)-\partial_x f(x,y)}{h}$ as the left $y$-derivative of $\partial_x f$ (if exists).

We can show that the updating rule  \eqref{eq:updated_policy} always improves in terms of the value function. We refer to  Section \ref{sec:AppD} for the proof.

\begin{theorem}[Policy improvement]\label{prop:policy_improvement}  Assume Assumptions \ref{ass:Pi} and \ref{ass:g} hold. 
In addition, assume that $\widetilde{g}$ updated according to \eqref{eq:updated_policy} is increasing on $[0,\hat{x}_{\widetilde{g}}]$. Then it holds that
    \begin{eqnarray}\label{eq:policy_improvement}
        V^\lambda_{\widetilde{g}}(x,y)\ge V^\lambda_g(x,y),
    \end{eqnarray}
    where $V^\lambda_g, V^\lambda_{\widetilde{g}}$ are value functions associated with policies $\xi^g, \xi^{\widetilde{g}}$, see \eqref{eq:parameterized_policy}.
\end{theorem}

Among the conditions listed in the theorem, the requirement that 
$\widetilde{g}$ updated via \eqref{eq:updated_policy} is increasing on 
$[0,\hat{x}_{\widetilde{g}}]$ should be viewed as a propagated property rather than an assumption: Theorem~\ref{thm:policy_convergence}
establishes that the update preserves this monotonicity, so it suffices to verify it at initialization.

\begin{figure}[H]
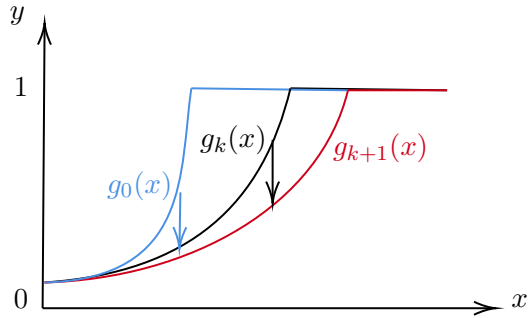

\centering
\include{figures/g_update_v3}
\caption{\label{figure:update} Demonstration of the Policy Iteration Algorithm.}
\end{figure}

\vspace{5pt}
Repeat the updating procedure \eqref{eq:updated_policy} iteratively, we have the general heuristic described in Algorithm \ref{alg1}. Although new in the literature, Algorithm \ref{alg1} can be broadly classified as a Policy Iteration method \cite{sutton2018reinforcement,thomas2016data}, which contains two steps: (1) a Policy Evaluation step that calculates the value function of a given policy ({\it cf.} line \ref{line:33} of Algorithm \ref{alg1}) and (2) a Policy Improvement step that aims to update the policy in the direction of improving the value function ({\it cf.} line \ref{line:44} of Algorithm \ref{alg1}).
The idea behind our specific design is that, we start with an initial $g_0(x)$ that has a sufficiently large exploration region $\mathcal{E}(g_0)$ ({\it i.e.}, Assumption \ref{ass:initial_policy}). Then, in each iteration $k$, we update  $g_k$ function  ``downwards'' into $\mathcal{E}(g_k)$ according to the Hessian information ({\it cf.} line \ref{line:33} of Algorithm \ref{alg1}), in which $u_k$ has a better regularity ({\it cf.} line \ref{line:44} of Algorithm \ref{alg1}). See Figure \ref{figure:update} for a demonstration.

\begin{algorithm}[h]
\caption{Policy Iteration for Entropy-regularized Optimal Stopping (PI-$\lambda$-OS)} \label{alg1}
\begin{algorithmic}[1]
\State Initialize $g_0(x)$ for $x\in[0,\infty)$ according to Assumption \ref{ass:initial_policy}.
\For{$k=0,1,\cdots,K-1$}
\State\label{line:33} Find $u_k(x,y) $  a ${C}^1(\mathbb{R}_+\times [0,1])\cap{C}^2\left(\overline{\mathcal{E}(g_k)}\right)$ solution to the following equations:
  \begin{eqnarray}
       &&(\mathcal{L}_x-\rho) u + \Big(\pi(x)-\rho \kappa \Big)y - \lambda y \log y                   =0\quad \text{on} \quad\mathcal{E}(g_k),\label{eq:uk-1}\\
     && -u_y = 0, \quad \text{on}\quad \mathcal{S}(g_k).\label{eq:uk-2}
  \end{eqnarray}
 \State \label{line:44} Update the strategy
 \begin{equation*}
      g_{k+1}(x) =
      \begin{cases}
       & \max \Big\{y< g_k(x)\,\Big| \partial_{xy} u_k(x,y) = 0 \Big\} \,\,\text{if} \,\, \partial_{xy}^- u_k(x,g_k(x))<0, \\
    &  g_{k+1}(x) =  g_k (x)  \qquad \text{otherwise.}  
      \end{cases}
 \end{equation*}
\EndFor
\end{algorithmic}
\end{algorithm}
 This iterative scheme is inspired by \cite{kumar2004numerical} for a one-dimensional singular control problem. However, the convergence result in that work requires a second-order condition to hold throughout the entire iteration process, which is difficult to verify.

In addition to the policy improvement result in Theorem \ref{prop:policy_improvement}, we provide the following policy convergence result under additional conditions on the initialization. 

\begin{assumption}[Initial policy]\label{ass:initial_policy}
Assume the initial policy $g_0$ satisfies the following conditions:
\begin{itemize}
   \item[{\bf (a)}] $g_0\in {C}^1([0,\hat{x}_{g_0}])$ is strictly increasing on $[0,\hat{x}_{g_0}]$,
    \item[{\bf (b)}] $g_0(0) = e^{-(1+
\frac{\kappa \rho}{\lambda})}$,  and
\item[{\bf (c)}] $- \alpha_- \Big(\kappa  +\frac{\lambda}{\rho}  \log (g_0(x))+ \frac{\lambda}{\rho}\Big) +\alpha_-H_\pi(x) - H_\pi^{\prime}(x)\cdot x \ge 0$ on $[0,\hat{x}_{g_0}]$.
\end{itemize}
\end{assumption}
Note that {\bf (c)} implies that $\hat{x}_{g_0} \leq \hat{x}_{g_\lambda}$. 

\begin{remark}[Justification of Assumption \ref{ass:initial_policy}] Assumption \ref{ass:initial_policy}-(a) is easy to satisfy and   Assumption \ref{ass:initial_policy}-(c) requires $\mathcal{E}(g_0)$ to be sufficiently large so that $\mathcal{E}(g_0) \supseteq \mathcal{E}(g_\lambda)$ . Assumption \ref{ass:initial_policy}-(b) seems to be the most restrictive one, as it requires the knowledge of $g_\lambda(0)=\exp(-(1+\frac{\kappa \rho}{\lambda}))$. In the model-free setting (such as in Section \ref{sec:model-free}) and when $g_\lambda(0)$ is not available, the following algorithm can be applied to learn $g_\lambda(0)$ at a fast convergence rate.

\begin{algorithm}[H]
\caption{Learning Initial Value $g_0(0)$ \label{alg0}}
\begin{algorithmic}[1]
\State Initialize $y_0 \in (0,1)$, $c_0>0$ and $\{\eta_i\}_{i \ge 1}$.
\For{$i=1,\cdots,$}
\State Acquire the value functions $u_i^+ = J^\lambda(0,y_{i}+\varepsilon_i;\xi=0)$ and $u_i^- =J^\lambda(0,y_{i}-\varepsilon_i;\xi=0)$,   with $\varepsilon_i = \min\{y_i,1-y_i, c_0/i\}$.
\State Update $y_i$ using the following zeroth order gradient descent with a two-point estimator:
\begin{eqnarray}
    y_{i+1} = y_i - \eta_i \frac{u_i^+-u_i^-}{2\varepsilon_i}. 
\end{eqnarray}
\EndFor
\end{algorithmic}
\end{algorithm}
We know that $J^\lambda(0,y;\xi=0) =  - \kappa y - \frac{\lambda y}{\rho} \log y$ is convex in $y$ with the minimizer taking value $y^*:=e^{-(1+\frac{\kappa \rho}{\lambda})}$. When the functional form (hence the minimizer) is unknown to the decision-maker, Algorithm \ref{alg0} converges linearly to the minimizer \cite{duchi2015optimal} when taking $\eta_i = \frac{\eta}{\sqrt{i}}$ with some $\eta>0$. Mathematically,
\begin{eqnarray}
    |y_i-y^*|^2\sim \mathcal{O}(\delta^i),
\end{eqnarray}
for some $\delta\in(0,1)$ depending on model parameters.

\end{remark}

\begin{remark}[Examples that satisfy Assumption \ref{ass:initial_policy}]\label{rmk:initial_policy}Assumption \ref{ass:initial_policy} is easy to satisfy and we provide two examples under the special case $\pi(x) = x^\theta$ with $0<\theta<1$: one is linear initialization and the other one is exponential initialization. For linear initialization, we can take 
\begin{eqnarray}
    g_0(x) = \min\left\{e^{-(1+\frac{\kappa\rho}{\lambda})} +2 (1-e^{-(1+\frac{\kappa\rho}{\lambda})})\Big(\frac{- \alpha_- \big(\kappa  +\frac{\lambda}{\rho}\big)}{\theta-\alpha_-}\Big)^{\theta}x, 1\right\},
\end{eqnarray}
which interpolates $(x_0,y_0)$ and $(x_1,y_1)$ with
$x_0=0$, $y_0=e^{-(1+\frac{\kappa\rho}{\lambda})}$, $x_1 =  \frac{1}{2}\Big(\frac{- \alpha_- \big(\kappa  +\frac{\lambda}{\rho}\big)}{\theta-\alpha_-}\Big)^{-\theta}$ and $y_1=1$. 

For exponential initialization, we can take 
\begin{eqnarray}
    g_0(x) = \min\left\{\exp\Big(\frac{\rho (\zeta-\alpha_-) (x)^{\zeta}}{- \alpha_-\lambda} - \frac{\rho}{\lambda}( \kappa+  \frac{\lambda}{\rho} )\Big), 1\right\},
\end{eqnarray}
for some $\zeta\in(\theta,1)$. This is inspired by the form of \eqref{eq:g_lam}, that under $\pi(x)=x^\theta$ gives
$$g_\lambda(x) =\min\left\{\exp\Big(\frac{\rho (\theta-\alpha_-) (x)^{\theta}}{- \alpha_-\lambda} - \frac{\rho}{\lambda}( \kappa+  \frac{\lambda}{\rho} )\Big), 1\right\}.$$
\end{remark}

\begin{theorem}[Policy convergence]\label{thm:policy_convergence} Let Assumptions \ref{ass:Pi}, \ref{ass:g} and \ref{ass:initial_policy} hold. Then we have the following results
\begin{eqnarray}
  \label{eq:policy_convergence1}  &&\lim_{k \rightarrow \infty} g_k = g_\lambda,\\
 \label{eq:policy_convergence2}   &&\lim_{k \rightarrow \infty} V^\lambda_{g_k} = V_{\lambda}.
\end{eqnarray} 
with $g_\lambda$ defined in \eqref{eq:g_lam}. In particular, we have $g_k$ is strictly increasing on $[0,\hat{x}_{g_k}]$ and $g_k\in C^1([0,\hat{x}_{g_k}])$ for all $k\in \mathbb{N}_+$.
\end{theorem}
Note that the monotonicity and smoothness of the boundary $g_k$ ensures the existence and uniqueness of the solution to \eqref{eq:uk-1}-\eqref{eq:uk-2}, which coincides with the value function under policy $\xi_{g_k}$ (see Theorem \ref{thm:property_g}). The proof of the previous theorem about policy convergence can be found in Section \ref{sec:AppE}.

\begin{remark} It is worth noting that
Theorem \ref{thm:policy_convergence} presents the first policy improvement result for RL algorithms related to free boundary problems. Existing literature primarily focuses on the development of continuous-time regular controls. Compared to policy improvement, showing policy convergence remains significantly more challenging, even for regular controls. To the best of our knowledge, the only results on policy convergence for regular controls are found in \cite{bai2023reinforcement, huang2022convergence, ma2024convergence, tran2024policy}, and all of these results are technically intricate.
\end{remark}

\subsection{Model-free implementation}\label{sec:model-free} 
When the model parameters are unknown, we are not able to calculate the solution $u_k$ of \eqref{eq:uk-1}-\eqref{eq:uk-2}.  In this case, we assume the learner has access to an environment simulator (see Algorithm \ref{alg3}). For each initial position $(x,y)$ and threshold function $g$ provided by the learner, the generator will return an instantaneous reward function associated with a random path $X^x$ and the corresponding control $Y^{y,\xi^g}$ (see line \ref{line:ins_rwd} in Algorithm \ref{alg3}). It is worth noting that the learner does not know the expression of the instantaneous reward function nor the generator of the dynamics.

\begin{algorithm}[H]
\caption{Simulator $\mathcal{G}$\label{alg3}}
\begin{algorithmic}[1]
\State {\bf Input:} Threshold function $g$, initial position $(x,y)$
\State {\bf Generate:}  Sample path $(X^{x},Y^{y,\xi^{g}})$ under policy $\xi^{g}$ (defined in \eqref{eq:parameterized_policy}) 
\State \label{line:ins_rwd} {\bf Return:} \begin{eqnarray}
 \qquad \qquad   \int_0^\infty e^{-\rho t} \Big(\Big( \pi(X^{x}_t)-\rho \kappa\Big)Y^{y,\xi^{g}}_t - \lambda Y^{y,\xi^{g}}_t \log (Y^{y,\xi^{g}}_t)\Big) \dd t\,\,
\end{eqnarray}
\end{algorithmic}
\end{algorithm}

\begin{algorithm}[H]
\caption{Sample-based Policy Iteration for Exploratory Optimal Stopping (SPI-$\lambda$-OS) }\label{alg2}
\begin{algorithmic}[1]
\State Initialize $g_0(x)$ according to Assumption \ref{ass:initial_policy}. Specify a grid size $\delta_x$ for partitioning the $x$-axis and a grid size $\delta_y$ for partitioning the $y$-axis. Also, specify an upper bound $\bar x:=N\delta_x$.
\For{$k=0,1,\cdots,K-1$}
\For{$x\in\{0,\delta_x,2\delta_x,\cdots,N\delta_x\}$ and $y\in \{0,\delta_y,2\delta_y,\cdots,\lfloor1/\delta_y\rfloor \}$}
\For{$m=1,\cdots,M$}
\State Acquire the simulator $u^{(m)}_k(x,y) = \mathcal{G}(x,y,g_k)$
\hfill \Comment{independent randomness across $M$ paths} 
 \EndFor
\State \label{line8} Calculate the approximated value function $\bar u_k(x,y) = \frac{\sum_{m=1}^M u^{(m)}_k(x,y)}{M}$
 
 \EndFor
\State \label{line10} Update the strategy
 \begin{equation*}
      g_{k+1}(x) =
      \begin{cases}
       & \max \Big\{y< g_k(x)\,\Big| \partial_{xy} \bar u_k(x,y) = 0 \Big\} \,\,\text{if} \,\, \partial_{xy}^- \bar u_k(x,g_k(x))<0, \\
    &  g_{k+1}(x) =  g_k (x)  \qquad \text{otherwise}  
      \end{cases}
 \end{equation*}
\EndFor
\end{algorithmic}
\end{algorithm}

By interacting with the generator, the learner approximates the value function by acquiring instantaneous rewards along multiple trajectories and take the average (see line \ref{line8} in Algorithm \ref{alg2}). We require independent randomness across $M$ paths. Mathematically, it means that the state processes are driven by independent Brownian motions.  The learner then implements a sample-based version of the Policy Improvement step (see line \ref{line10} in Algorithm \ref{alg2}). Note that the entire implementation avoids estimating model parameters. Instead, it iteratively updates the policy boundary. Therefore, this approach is referred to as a model-free implementation.


Notably, we do not directly estimate model parameters in Algorithm \ref{alg2}, which enhances robustness against model misspecification and environmental shifts \cite{agarwal2021theory,fazel2018global,hambly2021policy}.

\subsection{Numerical performance} We test the performance of Algorithm \ref{alg2} on a few examples.

For the underlying model, we set $\pi(x)=x^\theta$ with  $\theta=0.5$, $\mu=0.2$, $\sigma=0.2$, $\rho=0.5$ and $\kappa=5$ (assumed to be unknown to the learner). In the experiment, the learner set $\delta_x=\delta_y=0.02$, $\bar{x}=5$, and $M=20$.

\begin{figure}[H]
    \centering
    \includegraphics[width=7cm]{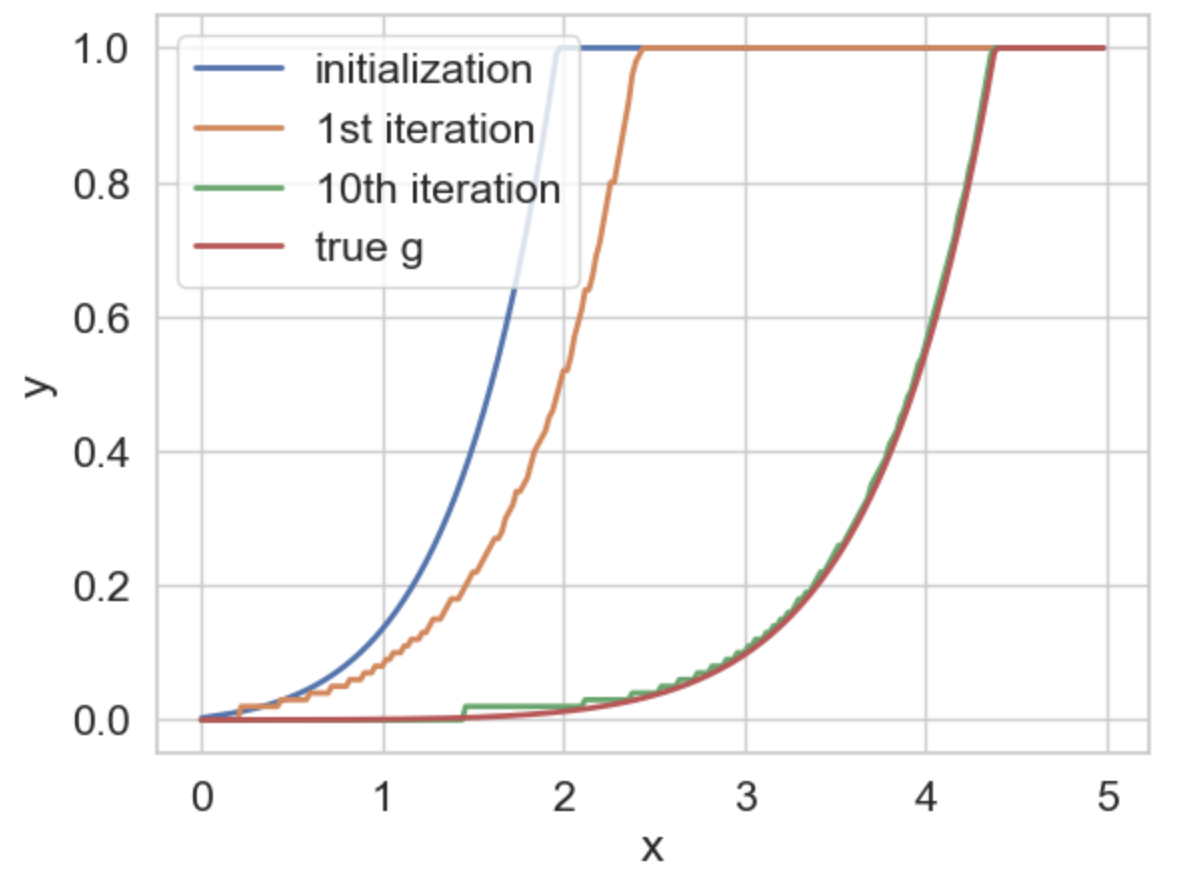}
\includegraphics[width=7cm]{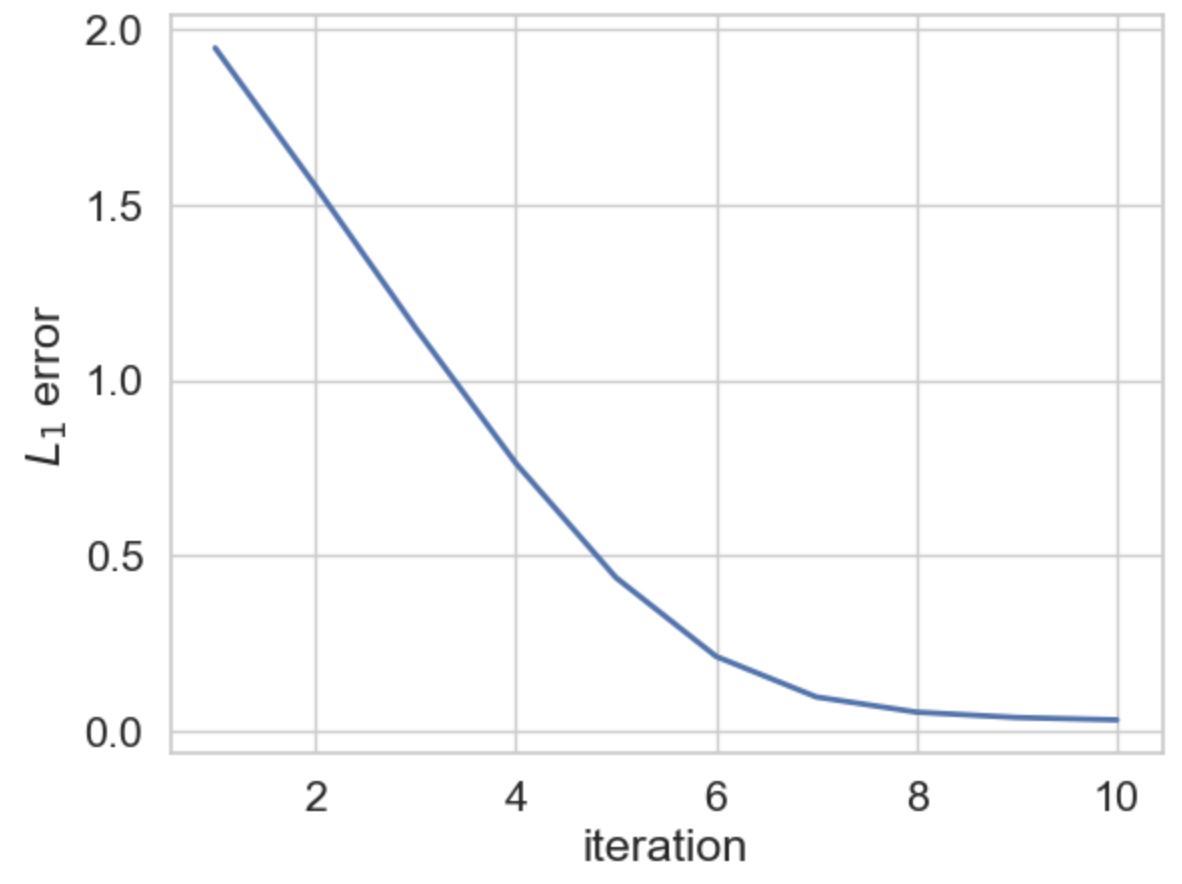}
    \caption{Exponential initialization. {\bf Left}: Ground-truth and learned $g$ function in selected iterations. {\bf Right}: Convergence to the ground truth in $L_1$ norm (outer iterations).}
    \label{fig:numerical1}
\end{figure}

\begin{figure}[H]
    \centering
    \includegraphics[width=7cm]{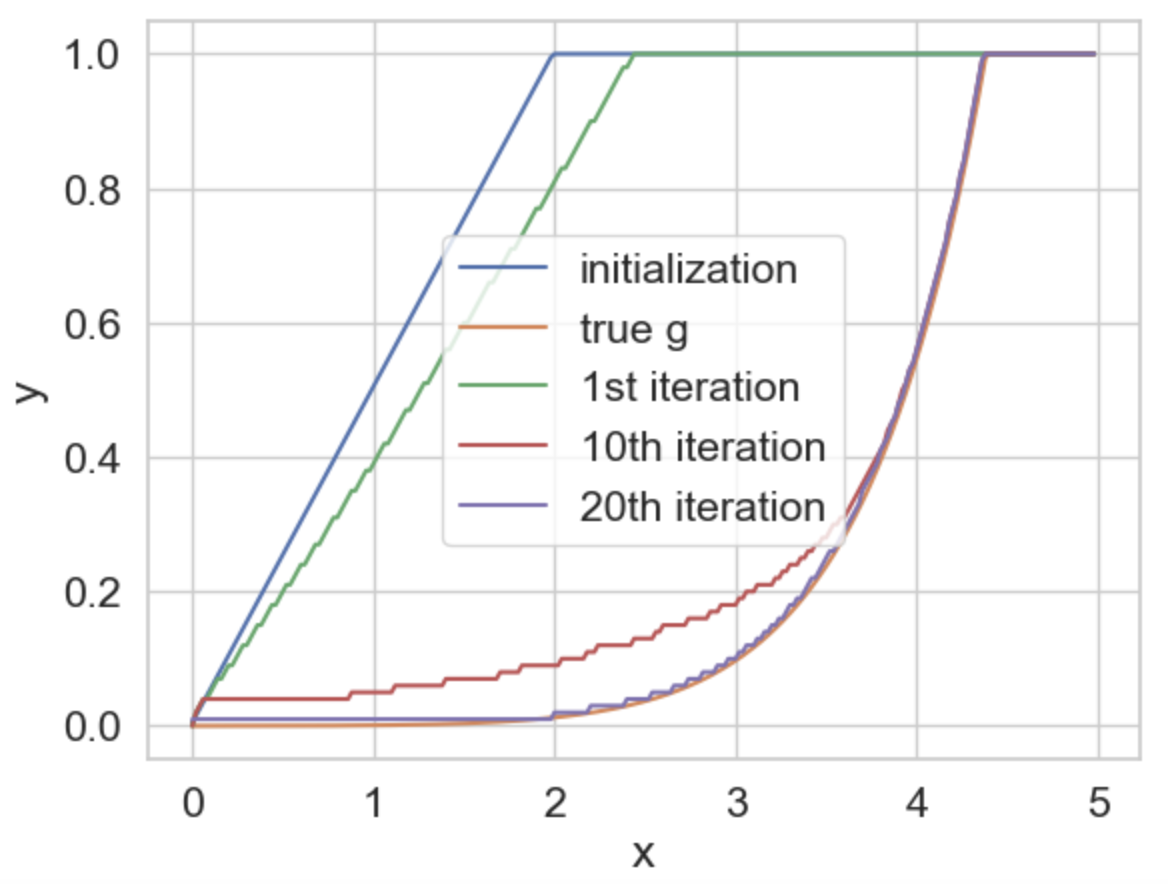}
\includegraphics[width=7cm]{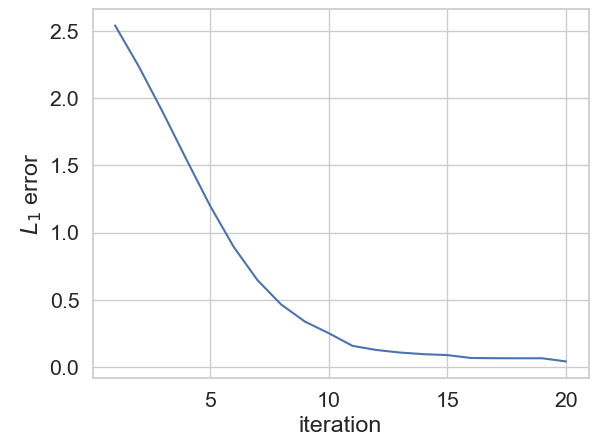}
    \caption{Linear initialization. {\bf Left}:  Ground-truth  and learned $g$ function in selected iterations . {\bf Right}: Convergence to the ground truth in $L_1$ norm (outer iterations).}
    \label{fig:numerical2}
\end{figure}

We evaluate the performance of the algorithm using two different initial policies: one with exponential initialization and the other with linear initialization (see Remark \ref{rmk:initial_policy}). As shown in Figure \ref{fig:numerical1}, Algorithm \ref{alg2} converges within $10$ outer iterations when using exponential initialization, highlighting its effectiveness. In addition, it takes approximately $20$ outer iterations for the algorithm to converge with linear initialization (see Figure \ref{fig:numerical2}). This suggests that learning the boundary is more challenging when $x$ is smaller.

\section{Detailed Proofs}
\label{Appendix A}

\subsection{Proof of Proposition \ref{prop xi does not explore}}
\label{sec:AppA}

For generic $\xi \in \mathcal A (1)$ and $z \in [0,1]$, we can define the stopping time $\tau ^\xi (z)$ as
$$
\tau^\xi (z) := \inf \{ t \geq 0 \, | \, \xi_t \geq z \},
$$
with the convention $\inf \emptyset := + \infty$.
Indeed, notice that $\tau^\xi (z) = + \infty$ if $z > \xi_\infty$.
By using Fubini's theorem and then the change of variable formula (see Chapter 0 in \cite{revuz.yor.2013continuous}), we have
$$
\begin{aligned}
    J^0(x&;\xi) \\
    =& \mathbb{E}\left[\int_{0}^{\infty} e^{-\rho t} \big(\pi (X^x_{t}) - \rho \kappa\big) \left(1 -\xi_\infty + \int_t^\infty  \dd \xi_s  \right) \dd t   \right] \\
     =& \mathbb{E}\left[ (1 -\xi_\infty) \int_{0}^{\infty} e^{-\rho t} \big(\pi (X^x_{t}) - \rho \kappa\big)   \dd t + \int_{0}^{\infty}  \left( \int_0^t e^{-\rho s} \big(\pi (X^x_{s}) - \rho \kappa\big) \dd s \right) \dd \xi_t \right] \\
     =& \mathbb{E}\left[ (1 -\xi_\infty) \int_{0}^{\infty} e^{-\rho t} \big(\pi (X^x_{t}) - \rho \kappa\big)   \dd t 
     +\int_{0}^{\xi_\infty} \bigg( \int_0^{\tau^\xi (z)} e^{-\rho s} \big(\pi (X^x_{s}) - \rho \kappa\big)   \dd s \bigg)   \dd z  \right] \\
     =& \int_{0}^1 \mathbb{E}\bigg[   \int_0^{\tau^\xi (z)} e^{-\rho s} \big(\pi (X^x_{t}) - \rho \kappa\big)   \dd s  \bigg]   \dd z 
     =\int_{0}^1 J(x;\tau^\xi (z) )  \dd z.
\end{aligned}
$$
Thus, since  $\tau^\xi(z)$ is an $\mathbb F ^W$-stopping time, by the optimality of $\tau^*$ we have
$$
J^0(x; \xi) \leq \int_{0}^1 \sup_{\tau \in \mathcal T } J(x;\tau )   \dd z = J(x;\tau^*) = V(x) = J^0(x;\xi^*), 
$$
where the last equality follows from the fact that  $\xi^*$ is such that $\tau^{\xi^*}(z) = \tau^*$ for any $z \in (0,1]$. This proves the first part of the proposition.

Notice that the optimal stopping time $\tau^*$ solving \eqref{eq real option OS} is unique. This is indeed given by $\tau^{*}=\inf\{t\geq0:\, X^{x}_t \leq b^*\}$, where the free boundary $b^*$ can be computed explicitly by employing the classical smooth-fit and smooth-pasting principles (see, e.g., \cite{DixitPindyck}).
Then, by repeating the previous argument with $\xi$ being optimal, we obtain that $J(x;\tau^\xi (z)) = J(x;\tau^*)$ $ \dd z$-a.e.\ in $(0,1)$. 
By uniqueness of the optimal stopping time, we deduce that $\tau^\xi (z) = \tau^*$ $ \dd z$-a.e.\ in $(0,1)$, which in turn implies $\xi = \xi^*$.


\subsection{Proof of Theorem \ref{theorem verification real option}}
\label{sec:AppB}

\emph{Step 1: Constructing a candidate optimal solution.}
We  will now search for a candidate solution $u: (0,\infty) \times [0,1] \to \R$ and a nondecreasing function $ g_{\lambda}: (0,\infty) \to  [0,1]$ such that 
\begin{equation}\label{eq split guess HJB}
\begin{cases}
\left( \mathcal{L}_{x}-\rho\right) u + (\pi(x)-\rho \kappa) y -\lambda y \log y = 0, &0< y \leq g_{\lambda}(x), \\ 
-u_{y}=0, &  y>g_{\lambda}(x),
\end{cases}   
\end{equation}
and we will later verify that $u =V^\lambda$ (see Theorem \ref{theorem verification real option} below).

If $u$ and $g_\lambda$ satisfy \eqref{eq split guess HJB},  for $0<y \leq g_\lambda(x)$,
 we can integrate the equation to obtain
\begin{equation}\label{eq:alpha_pm}
u(x, y)=
A_{2}(y) x^{\alpha_-}
+A_{1} (y)x^{\alpha_{+}}
+ \mathbb{E} \left[ \int_{0}^{\infty} e^{- \rho t} \big( (\pi(X^x_{t})-\rho \kappa ) y- \lambda y \log y \big) \dd t \right],
\end{equation}
with $\alpha_+, \alpha _-$ solving
$
\frac{1}{2} \sigma^{2} \alpha(\alpha-1)+\mu \alpha-\rho =0, 
$
and such that $\alpha_+ > 1 $ (as $\rho >\mu$) and $\alpha_-<0$.
Moreover, since $\rho > \mu$ and Assumption \ref{ass:Pi} hold, we  obtain that
\begin{equation*}\label{eq estimate V growth real option}
\begin{aligned}
V^\lambda (x,y) 
& \leq C  +  C 
 \mathbb E \bigg[ \int_0^{\infty} e^{-\rho t} (X^x_t)^\theta d t \bigg] 
 = C  + C x^\theta  
 \mathbb E \bigg[ \int_0^{\infty} e^{(-\rho + \theta (\mu + \frac{\sigma^2 }{2} (\theta-1) ))t} d t \bigg] 
 \leq C (1+x^\theta).
\end{aligned}
\end{equation*}
Thus, we take $A_{1} \equiv 0$ in \eqref{eq:alpha_pm}, as otherwise $u$ would be superlinear as $x \uparrow {\infty}$, not matching the growth behavior of $V^\lambda$.
Hence,  for $0<y \leq g_\lambda(x)$, a direct integration leads to
\begin{equation*}\label{constant_c}
u(x, y)
=A_{2}(y) x^{\alpha_-} + H_\pi(x)y - \kappa y -\frac{\lambda}{\rho} y \log y,\,\, \text{ with } \,\,H_\pi(x)  = \mathbb{E}\left[\int_0^\infty e^{-\rho t}\pi(X_t^x)\dd t\right].
\end{equation*}
To summarize, we have obtained the following candidate value function
$$
u(x, y) =
\begin{cases}
A_{2}(y) x^{\alpha_-}+ H_\pi(x)y -\kappa y-\frac{\lambda}{\rho} y \log y, & 0<y\leq g_{\lambda}(x), \\
A_{2}\left(g_{\lambda}(x)\right) x^{\alpha_-}+\pi(x)\,y-\kappa g_{\lambda}(x)-\frac{\lambda}{\rho} g_{\lambda}(x) \log \left(g_{\lambda}(x) \right), & y> g_{\lambda}(x).
\end{cases} 
$$ 
In order to determine $A_{2}$ and $g_{\lambda}$ we impose
\begin{align*}
 \lim_{z \uparrow g_{\lambda}(x) } u_{y}\left(x, z\right)=\lim_{z \uparrow g_{\lambda}(x) } u_{y x}\left(x, z\right)=0.
\end{align*}
From these conditions, we derive the system
\begin{equation*}
\begin{cases}
    A_{2}^{\prime}(g_{\lambda}(x))\,x^{\alpha_-} + H_\pi(x)-\kappa-\frac{\lambda}{\rho}(1+\log g_\lambda(x))=0, \\
    \alpha_{-} A_{2}^{\prime}(g_{\lambda}(x))x^{\alpha_- -1}+H'_\pi(x) = 0,
\end{cases}
\end{equation*}
which can be solved as
\begin{equation}\label{eq:optimal_g}
\begin{cases}
g_{\lambda}(x) = \exp \left( \frac{-\frac{H_\pi'(x)}{\alpha_-}x+H_\pi(x)-\kappa - \frac{\lambda}{\rho}}{\frac{\lambda}{\rho}}\right) \\
A_{2}(y) =\int_{g_\lambda(0)}^y  \frac{\kappa  +\frac{\lambda}{\rho}  \log (u)+ \frac{\lambda}{\rho}- H_\pi(g_\lambda^{-1}(u))}{(g_\lambda^{-1}(u))^{\alpha_-}} \dd u.
\end{cases}
\end{equation} 
Notice that 
\begin{eqnarray}
   g_{\lambda}'(x) =  g_{\lambda}(x)  \left(\frac{-1}{\alpha_-} (H_\pi'(x)x)'+H'_\pi(x)\right) > 0,
\end{eqnarray}
as $(H_\pi'(x)x)' > 0$ and $H'_\pi(x) \ge 0$ according to Assumption \ref{ass:Pi}.
Hence $g_\lambda$ is strictly increasing on $(0,\hat{x}_{g_\lambda}]$ with $\hat{x}_{g_\lambda} = \min\{x\in \mathbb{R}_+: g_\lambda(x)=1\}$. 

Notice that $g_\lambda$ delimits two regions only for $y \geq y^{\lambda}$, with
\begin{equation}\label{eq definition of minimal y}
 y^\lambda:= g_\lambda(0) = e^{-1-\frac{\kappa \rho}{\lambda} }  \in(0,1).
\end{equation}

For $y<y^\lambda$ we notice that the function
$$
u_0 (x, y)
:=  H_\pi(x)y - \kappa y - \frac{\lambda}{\rho} y \log y,
$$ 
is such that 
\begin{equation}\label{eq def candidate below threshold}
\begin{cases}
    (\mathcal{L}_x -\rho) u_{0}+(x^\theta-\rho \kappa) y-\lambda y \log y=0, &  \\ 
-     \partial_{y} u_{0}= -H_\pi(x)+\left( \kappa +\frac{\lambda}{\rho}(\log y+1)\right)<0, & \text { for } y<y^{\lambda}.
\end{cases}
\end{equation}
Thus, 
$u_{0}(x, y)$ solves the HJB and for $y<y^{\lambda}$ and we expect $\xi^{\lambda} \equiv 0$ to be optimal.

Summarizing the previous heuristic,  
set
$ F(x, y):=A_{2} (y) \,x^{\alpha_-} + H_\pi(x)y - \kappa y -\frac{\lambda}{\rho} y \log y$ and define the candidate value as
\begin{equation*}
u(x,y) :=
\begin{cases}
H_\pi(x)y - \kappa y - \frac{\lambda}{\rho} y \log y, & y<y^\lambda, \\
F(x, y), & y\leq g_\lambda (x), \ y \geq y^\lambda, \\
F(x, b_\lambda^{-1} (x) ), & y> g_\lambda (x), \ y \geq y^\lambda, 
\end{cases}
\end{equation*}
and we extend $g_\lambda$ on $[0,\infty)$ as follows
\begin{equation*}
g_{\lambda}(x)= 
\begin{cases}
g_{\lambda}(x) ,  &  x \in [0, \hat{x}_{g_\lambda}] \\ 
1, & x \ge \hat{x}_{g_\lambda}.
\end{cases}
\end{equation*}
\vspace{0.25cm}

\emph{Step 2: Regularity of the candidate value function.}
In this step we show the $C^2$-regularity of $u$.
We split $(0,\infty) \times [0,1]$ into the three sets $\mathcal E  _0 := \{ (x, y) \, |\, y<y^\lambda \}$, $\mathcal E := \{ (x, y) \, |\, y \geq y^\lambda, \ y<g_{\lambda} (x)\}$ and $ \mathcal S := \{ (x, y) \, |\, y \geqslant g_{\lambda}(x), \ y \geqslant y^\lambda\}$. 
Clearly, $u$ is regular in the interior of $\mathcal E _0$, $\mathcal E$ and $\mathcal S$.
In order to show the $C^2$-regularity across $\mathcal E _0$ and $\mathcal E$, it is sufficient to notice that $A_2(y^\lambda)=A_2'(y^\lambda)=A_2''(y^\lambda)=0$.

We next study the behavior across $ \mathcal S$ and $\mathcal E$. 
To this end, recall the definition of 
$ F$ and notice that
$$
\begin{cases}
u(x, y)=F(x, y), & y<g_{\lambda}(x), \\ 
u(x, y)=F\left(x, g_{\lambda}(x)\right), & y \geqslant g_{\lambda}(x).
\end{cases}
$$
Clearly we have
$$
\lim _{y \uparrow g_{\lambda}(x) }  u (x, y)  = F \left(x, g_{\lambda}(x)\right)  =\lim _{y \downarrow g_{\lambda}(x)} u(x, y),
$$
so that $u$ is continuous across the boundary.
Furthermore, after straightforward calculations we find
$$
\begin{aligned}
& u_{x}(x, y)=F_{x}(x, y), \quad y<g_{\lambda}(x), \\
& u_{x}(x, y)= F_{x} \left( x, g_{\lambda}(x) \right) + \underbrace{F_{y}\left(x, g_{\lambda}(x)\right)\left(g_{\lambda} \right)^{\prime}(x)}_{=0}, \quad y>g_{\lambda}(y),
\end{aligned}
$$
so that
$$
\lim _{y \uparrow g_{\lambda}(x)} u_{x}(x, y) =F_{x}\left(x, g_{\lambda}(x)\right)=  \lim _{y \downarrow g_{\lambda}(x)} u _x (x, y), 
$$
which  in turn implies the continuity of $u_x$.
Finally, 
\begin{align*}
& u_{x x}(x, y)=F_{x x}(x, y), \quad y<g_{\lambda}(x), \\
& u_{x x}(x, y)=F_{x x}\left(x, g_{\lambda}(x)\right)+\underbrace{F_{x y}\left(x, g_{\lambda}(x)\right)\left(g_{\lambda}\right)^{\prime}(x)}_{=0}, \quad y>g_{\lambda}(x), 
\end{align*}
from which we conclude that
\begin{equation} 
 \lim _{y \downarrow g_{\lambda}(x)} u_{x x}(x, y)=F_{x x}\left(x, g_{\lambda}(x)\right)=\lim _{y \uparrow g_{\lambda}(x)} u_{x x}(x, y),
 \end{equation}
which is the continuity of $u_{xx}$.
For the derivatives involving $y$, by using the definitions of $A_2$ and $g_\lambda$ in \eqref{eq:optimal_g}, a direct computation leads to
$$
\begin{aligned}
& \lim _{ x \downarrow g^{-1}_{\lambda}(y) } u_{y}(x, y)
=F_{y}\left( g^{-1}_{\lambda}(y), y\right) = 0 = \partial _y F\left(x, g_{\lambda}(x)\right) , \\
& \lim _{ x \downarrow g^{-1}_{\lambda}(y) } u_{yy}(x, y)
=F_{yy}\left(g^{-1}_{\lambda}(y), y\right) = 0 = \partial _{yy}^2 F\left(x,g_{\lambda}(x)\right) , \\
 & \lim _{ x \downarrow g^{-1}_{\lambda}(y) } u_{yx}(x, y)
=F_{yx}\left( g^{-1}_{\lambda}(y), y\right) = 0 = \partial _{yx}^2 F\left(x, g_{\lambda}(x)\right) ,
\end{aligned}
$$
which gives us the desired regularity.
\vspace{0.25cm}

\emph{Step 3: The candidate value function solves the HJB equation.} In this step we show that $u$ solves the HJB equation. Let $\mathcal E_0$, $\mathcal E $ and $\mathcal S$ be as in the previous step.
First, on the set $\mathcal E  _0$,  $u$ satisfies the HJB equation by  \eqref{eq def candidate below threshold}. 
Secondly, on the set $\mathcal E $, by construction we have
$(\mathcal{L}_x- \rho) u +( \pi(x)- \rho \kappa) y- \lambda y \log y                  =0,$
and we need to verify that $-u_y \leq 0$. 
To this end, use the definition of $A_2$ and $b_\lambda$ to compute the derivative 
\begin{align*}
  - u_{yx} (x,y) &=   -\alpha_-A_{2}' (y) x^{\alpha_--1} - H'_\pi(x)\\
  &=-\alpha_-\frac{\kappa  +\frac{\lambda}{\rho}  \log (y)+ \frac{\lambda}{\rho}- H_\pi(g_\lambda^{-1}(y))}{(g_\lambda^{-1}(y))^{\alpha_-}}  x^{\alpha_--1} - H'_\pi(x)\\
   &=-\alpha_-\frac{\big(\kappa  +\frac{\lambda}{\rho}  \log (y)+ \frac{\lambda}{\rho}- H_\pi(g_\lambda^{-1}(y))\big) x^{\alpha_--1}+\frac{1}{\alpha_-}H'_\pi(x)(g_\lambda^{-1}(y))^{\alpha_-}}{(g_\lambda^{-1}(y))^{\alpha_-}}  x^{\alpha_--1} \\
   &= -\alpha_-\frac{\frac{H_\pi'(g_\lambda^{-1}(y))}{-\alpha_-} g_\lambda^{-1}(y)g_\lambda^{-1}(y)\Big(x^{\alpha_--1}-(g_\lambda^{-1}(y))^{\alpha_--1}\Big)}{(g_\lambda^{-1}(y))^{\alpha_-}}  x^{\alpha_--1} \leq 0
\end{align*}
where the last inequality holds since $x^{\alpha_--1}-(g_\lambda^{-1}(y))^{\alpha_--1}< 0$, which is a direct consequence of $g_\lambda^{-1}(y)<x$.
Thus, since $-u_y(x,g_\lambda(x)) = 0$, we have that $-u_y \leq 0$ on $\mathcal E$.

We next study the behavior  on $ \mathcal S$. 
By construction we have $-u_y=0$ in $\mathcal S$.
Moreover, using the expression for $u, u_{x}$ and $u_{xx}$ computed in the previous step, for any $x$ such that  $y>g_{\lambda}^{-1}(x)$ we find
$$
\begin{aligned}
(\mathcal{L}_{x}-\rho) u & + (\pi(x)-\rho \kappa) y -\lambda y \log y \\
&=  \frac{1}{2} \sigma^{2} x^{2} \partial_{x x}\left(F\left(x, b_{\lambda}^{-1}(x)\right)\right) + \mu x \partial_{x}\left(F\left(x, b_{\lambda}^{-1}(x)\right)\right) 
-\rho F\left(x, b_{\lambda}^{-1}(x)\right) \\
& \quad+(\pi(x)-\rho \kappa) y -\lambda y \log y                   \\
& =\frac{1}{2} \sigma^{2} x^{2} F_{x x}\left(x, b_{\lambda}^{-1}(x)\right)+\mu x F_{x}\left(x, b^{-1}(x)\right)-\rho F\left(x, b^{-1}(x)\right) \\
& \quad+(\pi(x)-\rho \kappa) y -\lambda y \log y.
\end{aligned}
$$
Hence, using the ordinary differential equation for $y<g_{\lambda}(x)$, in the limit as $y \uparrow g_{\lambda}(x)$, we proceed with
$$ 
\begin{aligned}
(\mathcal{L}_{x}-\rho) u & + (\pi(x)-\rho \kappa) y -\lambda y \log y \\
& =-(\pi(x)-\rho \kappa) g_{\lambda}(x)+\lambda g_{\lambda}(x) \log \left(g_{\lambda}(x)\right)+(\pi(x)-\rho \kappa) y-\lambda y \log y                  \\
& =\int_{g_{\lambda}(x)}^{y}\left( \pi(x)-\rho \left(k+\frac{\lambda}{\rho}(1+\log z)\right)\right) d z \\
& <\left(y-g_{\lambda}(x)\right)\left(\pi(x)-\rho\left(k+\frac{\lambda}{\rho}\left(1+\log \left(g_{\lambda}(x)\right)\right)\right)\right) \\
& =\underbrace{\left(y-g_{\lambda}(x)\right)}_{ \geq 0} \left( \pi(x) -\rho\left(\frac{H_\pi'(x)}{-\alpha}x +H_\pi(x) \right)\right).
\end{aligned}
$$
Note that $H_\pi$ satisfies the following ODE:
\begin{eqnarray}
    \mu x\, H'_\pi(x) -\rho H_\pi(x) +\pi(x) = -\frac{1}{2}\sigma^2 x^2 H^{\prime \prime}_\pi(x).
\end{eqnarray}
Therefore
\begin{eqnarray}
 \pi(x) -\rho\left(\frac{H_\pi'(x)}{-\alpha}x +H_\pi(x) \right)= - \left(\mu + \frac{\rho}{-\alpha_-}\right) x\, H'_\pi(x)-\frac{1}{2}\sigma^2 x^2 H^{\prime \prime}_\pi(x) \leq 0,
\end{eqnarray}
where the last inequality holds since  $\frac{1}{2} \sigma^{2} \alpha_{-}\left(\alpha_{-}-1\right)+\mu \alpha_{-}-\rho=0$ and $\alpha_-<0$. Hence $(\mathcal{L}_{x}-\rho) u  + (\pi(x)-\rho \kappa) y -\lambda y \log y\leq 0$ on $\mathcal{S}$, so that $u$ solves the HJB equation in the region $\mathcal S$.
\vspace{0.25cm}

\emph{Step 4: Verifying the actual optimality of the candidate solution.} The last step consists of constructing a candidate optimal control $\xi^\lambda$ as in the statement of the theorem and verifying its optimality and that of the candidate value function $u$ by the help of It\^o-Meyer's formula. The argument is nowadays classical (see, e.g., Chapter VIII in\cite{fleming.soner2006} or \cite{angelis&ferrari&moriarty2019} for a related setting), so that we omit the details.

\subsection{Proof of Theorem \ref{thm:property_g}}
\label{sec:AppC}
{Note that
     $g^{-1}$ is well defined under Assumption \ref{ass:g} with $g^{-1}(x)>0$ for all $x\in \mathbb{R}_+$.  We seek for a semi-explicit solution to
    \begin{eqnarray}
      &&(\mathcal{L}_x-\rho) u + \Big(\pi(x)-\rho \kappa \Big)y - \lambda y \log y                   =0\quad \text{on} \quad\mathcal{E}(g),\label{eq:HJB_g_1}\\
     && -u_y = 0, \quad \text{on}\quad \mathcal{S}(g).\label{eq:HJB_g_2}
  \end{eqnarray}
 Similar to the analysis in Section \ref{sec:real_option} and following the idea in \cite{ferrari2015integral,ferrari2016irreversible}, take
  \begin{eqnarray}\label{eq:u_g}
      u(x,y) = A(y) x^{\alpha_-} + H_\pi(x)\,y - \kappa y -\frac{\lambda}{\rho} y \log y                  ,
  \end{eqnarray}
  with $\alpha_-$ defined in \eqref{eq:alpha_pm} and
$H_\pi(x)$ defined as 
\begin{eqnarray} \label{eq:R_Pi}
H_\pi(x) := \mathbb{E}\Big[\int_{0}^\infty e^{-\rho t}\pi(X^x_t)\dd t\Big].
\end{eqnarray}
Under Assumption \ref{ass:Pi}, we have $H_\pi(x)\in {C}^2(\mathbb{R}_+)$ with sublinear growth \cite[Eqn (2.11)]{menaldi1983some}.

Then,
  \begin{eqnarray}
      u_y(x,y) = A'(y) x^{\alpha_-} +H_\pi(x)  - \kappa  -\frac{\lambda}{\rho}  \log y                  -\frac{\lambda}{\rho}. 
  \end{eqnarray}
 Setting $u_y(g^{-1}(y),y)=0$, we have
 \begin{eqnarray}
     A'(y) (g^{-1}(y))^{\alpha_-} + H_\pi(g^{-1}(y))  - \kappa  -\frac{\lambda}{\rho}  \log y                  -\frac{\lambda}{\rho} = 0.
 \end{eqnarray}
 Hence, $ A'(y) = \frac{\kappa  +\frac{\lambda}{\rho}  \log y                  + \frac{\lambda}{\rho}- H_\pi(g^{-1}(y))}{(g^{-1}(y))^{\alpha_-}}$ and
 \begin{eqnarray}\label{eq:A}
     A(y) = \int_{g(0)}^y  \frac{\kappa  +\frac{\lambda}{\rho}  \log (u)+ \frac{\lambda}{\rho}- H_\pi(g^{-1}(u))}{(g^{-1}(u))^{\alpha_-}} \dd u.
 \end{eqnarray}
It is easy to check that $u(x,y)$ is ${C}^1$ in $y$ on $\mathbb{R}_+\times[0,1]$ and that, under condition $g\in {C}^1([0,\hat{x}_g])$,  we have $u(x,y)\in {C}^{2}(\overline{\mathcal{E}(g)})$.

Now, on the curve $y=g(x)$,
\begin{eqnarray}
    u^+_x(x,y) = A(g(x)) \alpha_-x^{\alpha_--1} + R'_\pi(x) g(x). 
\end{eqnarray}
On the other hand,  we have $u(x,y) = u(x,g(x))$ for $(x,y)\in \mathcal{S}(g)$. Hence
\begin{eqnarray*}
   && \frac{u\big(x,g(x)\big)-u\big(x-\delta,g(x)
    \big)}{\delta} =  \frac{u\big(x,g(x)\big)-u\Big(x-\delta,g(x-\delta)\Big)}{\delta}\\
    &=& \frac{1}{\delta}\Big(A(g(x)) x^{\alpha_-}-A(g(x-\delta)) ((x-\delta))^{\alpha_-} \Big) + \frac{1}{\delta} \Big(H_\pi(x) g(x)-H_\pi(x-\delta) g(x-\delta)\Big)\\
    && - \frac{\kappa}{\delta}\Big( g(x)- g(x-\delta)\Big) -\frac{\lambda}{\rho\delta} \Big( g(x) \log \big(g(x)
    \big) - g(x-\delta) \log \big(g(x-\delta)
    \big)  \Big),
\end{eqnarray*}
and
\begin{eqnarray*}
   && u_x^-(x,y) = \lim_{\delta\downarrow 0} \frac{u\big(x,g(x)\big)-u\big(x-\delta,g(x)
    \big)}{\delta} \\
    &&= A(g(x)) \alpha_-x^{\alpha_--1} + R'_\pi(x)g(x) + u_y(x,g(x)) \\
    &&= A(g(x)) \alpha_-x^{\alpha_--1} + R'_\pi(x)g(x) g(x) = u_x^+(x,y),
\end{eqnarray*}
where the last equality holds since $u_y(x,g(x)) = 0$. Hence $u\in {C}^{1}(\mathbb{R}_+\times[0,1])$.

It should be noted that $u_{xx}$ fails to be continuous across the boundary
although it remains bounded on any compact subset of $\mathbb{R}_+\times[0,1]$.
}

Recall $Y_t^{ y,\xi^{g}}$ is the process  under the control $\xi^{ {g}}$,  defined in \eqref{eq:parameterized_policy}. 
Fix the initial condition $(X_0,Y_0^{ y,\xi^{g}})=(x,y)$ and take $R>0$. Set $\tau_R:=\inf\{t \ge 0: X_t \ge R\}$.  Applying Ito's formula in the weak version to $u$ (see \cite[Chapter 8, Section VIII.4, Theorem 4.1]{fleming2006controlled} and \cite[Theorem 4.2]{angelis2019solvable}),  up to the
stopping time $\tau_R\wedge T$, for some $T > 0$, to obtain
    \begin{eqnarray*}
     \mathbb{E}\Big[e^{-\rho (\tau_R\wedge T)} u(X^x_{\tau_R\wedge T},Y^{ y,\xi^g}_{\tau_R\wedge T})\Big] &=& u(x,y)  
     +\mathbb{E}\Big[\int_0^{\tau_R\wedge T} e^{-\rho t} \Big(\mathcal{L}u(X^x_t,Y^{ y,\xi^{g}}_{t-})-\rho V^\lambda_g(X^x_t,Y^{ y,\xi^{g}}_{t-})\Big) \dd t\Big] \\
      && + \mathbb{E}\Big[\int_0^{\tau_R\wedge T} e^{-\rho t} \, \partial_y u(X^x_t,Y^{ y,\xi^{g}}_{t-}) \dd Y_t^{ y,\xi^{g}} \Big]\\
      &&+\mathbb{E}\Big[\sum_{0\leq t\leq \tau_R\wedge T} \Big( \Delta u(X^x_t,Y^{ y,\xi^{g}}_{t}) - \partial_y u(X_t,Y^{ y,\xi^{g}}_{t-})\Delta Y_t^{ y,\xi^g} \Big)\Big],
  \end{eqnarray*}
  with the notations $\Delta u(X^x_t,Y^{ y,\xi^{g}}_{t}) =u(X^x_t,Y^{ y,\xi^{g}}_{t})-u(X^x_t,Y^{ y,\xi^{g}}_{t-})$ and $\Delta Y_t^{ y,\xi^{g}} =Y_{t}^{ y,\xi^{g}}-Y_{t-}^{ y,\xi^{g}} $.

Recall that  any admissible control $\xi^g$  can be decomposed into the sum of its continuous part and its pure jump part, i.e., $\dd \xi^g = \dd (\xi^g)^{\text cont}+\Delta \xi^g$. Hence we have the decomposition $\dd Y^{ y,\xi^g} = \dd (Y^{ y,\xi^g})^{\text cont}+\Delta Y^{ y,\xi^g}$, and therefore, 
\begin{eqnarray}
     \mathbb{E} \Big[e^{-\rho (\tau_R\wedge T)} u(X^x_s,Y^{y,\xi^{g}}_s)\Big] &=&  u(x,y)  
     +\mathbb{E}\Big[\int_0^{\tau_R\wedge T} e^{-\rho t} \Big(\mathcal{L}_x\, u(X^x_t,Y^{ y,\xi^g}_{t-})-\rho V^\lambda_g(X^x_t,Y^{ y,\xi^g}_{t-})\Big) \dd t \Big]\nonumber\\
      && +\mathbb{E}\Big[ \int_0^{\tau_R\wedge T} e^{-\rho t} \, \partial_y u(X^x_t,Y^{ y,\xi^g}_{t-}) \dd (Y^{ y,\xi^g})^{\text cont}_t\Big] \nonumber\\
      &&+\mathbb{E}\Big[\sum_{0\leq t\leq \tau_R\wedge T} \Delta u(X^x_t,Y^{ y,\xi^g}_{t}))\Big]\\
      &= &u(x,y) -\mathbb{E}\Big[\int_0^{\tau_R\wedge T} e^{-\rho t} \Big(\Big(\pi(X^x_t)-\rho \kappa\Big)Y_t^{ y,\xi^g} - \lambda  Y_t^{ y,\xi^g} \log\Big(Y_t^{ y,\xi^g}\Big) \Big)\dd t\Big].\label{eq:ito2} 
\end{eqnarray}
The last equation holds due to the following facts:
\begin{itemize}
    \item[(i)] It holds that 
    \begin{eqnarray*}
       && \mathbb{E}\Big[\int_0^{\tau_R\wedge T} e^{-\rho t} \Big(\mathcal{L}_x\,u(X^x_t,Y^{ y,\xi^g}_{t-})-\rho V^\lambda_g(X^x_t,Y^{ y,\xi^g}_{t-})\Big) \dd t \Big] \\&=& -\mathbb{E}\Big[\int_0^{\tau_R\wedge T} e^{-\rho t} \Big(\Big(\pi(X^x_t)-\rho \kappa\Big)Y_t^{ y,\xi^g} - \lambda  Y_t^{ y,\xi^g} \log\Big(Y_t^{ y,\xi^g}\Big) \Big)\dd t\Big],
    \end{eqnarray*}
    as pure jump could only possibly happen at time $0$,  $\xi^g$ keeps $(X^x_t,Y_t^{ y,\xi^g})$ within $\overline{\mathcal{E}(g)}$ for $t >0$ and the fact \eqref{eq:HJB_g_1};
    \item[(ii)] $\mathbb{E}\Big[\sum_{0\leq t\leq \tau_R\wedge T} \Delta u(X^x_t,Y^{ y,\xi^g}_{t}))\Big] = 0$ as $-u_y = 0$ on $\overline{\mathcal{S}(g)}$ according to the facts that \eqref{eq:HJB_g_2} holds
    and  $u\in{C}^{1}(\mathbb{R}_+\times[0,1])$; 
    \item[(iii)] $\mathbb{E}\Big[ \int_0^{\tau_R\wedge T} e^{-\rho t} \, \partial_y u(X^x_t,Y^{ y,\xi^g}_{t-}) \dd (Y^{ y,\xi^g})^{\text cont}_t\Big]=0$ as the  control increases only at the boundary $\overline{\mathcal{S}(g)}\cap \overline{\mathcal{E}(g)}$, where $-u_y = 0$.
\end{itemize}
When taking limits as $R \rightarrow\infty$, we have $\tau_R\wedge T\rightarrow T$, $\mathbb{P}$-a.s.  By standard properties of Brownian motion, it is easy to prove that the integral terms in the last expression on the right-hand side of \eqref{eq:ito2} are uniformly bounded in $L^2(\Omega, \mathbb{P})$, hence uniformly integrable. 
Moreover,
$u$ in \eqref{eq:u_g} has sublinear growth by straightforward calculation. Then we also take limits as $T\rightarrow \infty$ and it follows
that
\begin{eqnarray*}
    u(x,y) =\mathbb{E}\Big[\int_0^{\infty} e^{-\rho t} \Big(\Big(\pi(X^x_t)-\rho \kappa\Big)Y_t^{ y,\xi^g} - \lambda  Y_t^{ y,\xi^g} \log\Big(Y_t^{ y,\xi^g}\Big) \Big)\dd t\Big].
\end{eqnarray*}
Therefore, we have $u\equiv V^\lambda_g$, which completes the proof.


\subsection{Proof of Theorem \ref{prop:policy_improvement}}
\label{sec:AppD}

   First of all, $V^\lambda_g \in {C}^{2}\left(\overline{\mathcal{E}(g)}\right)$ according to Theorem \ref{thm:property_g}. Therefore, by the design of the algorithm, $\widetilde{g}$ satisfies the conditions in Assumption \ref{ass:g}. Hence Theorem \ref{thm:property_g} also applies to $\widetilde{g}$.

    Denote $Y_t^{y, \xi^{\widetilde{g}}}$ as the process under control $\xi^{ \widetilde{g}}$, with $ \widetilde{g}$ defined in  \eqref{eq:updated_policy}. Fix the initial condition $(x,y)$ and take $R>0$. Set $\tau_R:=\inf\{t \ge 0: X^x_t \ge R\}$. Similar to the proof of Theorem \ref{thm:property_g}, apply Ito's formula in the weak version to $V^\lambda_g$ (see \cite[Chapter 8, Section VIII.4, Theorem 4.1]{fleming2006controlled} and \cite[Theorem 4.2]{angelis2019solvable}) up to the
stopping time $\tau_R\wedge T$ for some $T > 0$. We then obtain
   {\allowdisplaybreaks \begin{eqnarray*}
     &&\mathbb{E}\Big[e^{-\rho (\tau_R\wedge T)} V^\lambda_g(X^x_{\tau_R\wedge T},Y^{y, \xi^{\widetilde{g}}}_{\tau_R\wedge T})\Big] \\
     &=& V^\lambda_g(x,y)  
     +\mathbb{E}\Big[\int_0^{\tau_R\wedge T} e^{-\rho t} \Big(\mathcal{L}_x V^\lambda_g(X^x_t,Y^{y, \xi^{\widetilde{g}}}_{t-})-\rho V^\lambda_g(X^x_t,Y^{y, \xi^{\widetilde{g}}}_{t-})\Big) \dd t\Big] \\
      && + \mathbb{E}\Big[\int_0^{\tau_R\wedge T} e^{-\rho t} \, \partial_y V^\lambda_g(X^x_t,Y^{y, \xi^{\widetilde{g}}}_{t-}) \dd Y_t^{ \widetilde{g}} \Big]\\
      &&+\mathbb{E}\Big[\sum_{0\leq t\leq \tau_R\wedge T}  e^{-\rho t}\Big( \Delta V^\lambda_g(X^x_t,Y^{y, \xi^{\widetilde{g}}}_{t}) - \partial_y V^\lambda_g(X^x_t,Y^{y, \xi^{\widetilde{g}}}_{t-})\Delta Y_t^{y, \xi^{\widetilde{g}}} \Big)\Big]\\
      &=& V^\lambda_g(x,y) -\mathbb{E}\Big[\int_0^{\tau_R\wedge T} e^{-\rho t} \Big(\Big(\pi(X^x_t)-\rho \kappa\Big)Y_t^{y, \xi^{\widetilde{g}}} - \lambda  Y_t^{y, \xi^{\widetilde{g}}} \log\Big(Y_t^{y, \xi^{\widetilde{g}}}\Big) \Big) \dd t\Big] 
      \\
      && +\mathbb{E}\Big[ \int_0^{\tau_R\wedge T} e^{-\rho t} \, \partial_y V^\lambda_g(X^x_t,Y^{y, \xi^{\widetilde{g}}}_{t-}) \dd (Y^{y, \xi^{\widetilde{g}}})_t^{\rm cont} \Big]\\
      && +\mathbb{E}\Big[\sum_{0\leq t\leq {\tau_R\wedge T}} e^{-\rho t}\Big( \Delta V^\lambda_g(X^x_t,Y^{y, \xi^{\widetilde{g}}}_{t-})  \Big)\Big]\\
      &\geq &V^\lambda_g(x,y) -\mathbb{E}\Big[\int_0^{\tau_R\wedge T} e^{-\rho t} \Big(\Big(\pi(X^x_t)-\rho \kappa\Big)Y_t^{y,\xi^{\widetilde{g}}} - \lambda  Y_t^{y, \xi^{\widetilde{g}}} \log\Big(Y_t^{y, \xi^{\widetilde{g}}}\Big) \Big)\dd t\Big],
\end{eqnarray*}}
\noindent with the notations $\Delta V^\lambda_g(X^x_t,Y^{y,\xi^{\widetilde{g}}}_{t}) =V^\lambda_g(X^x_t,Y^{y,\xi^{\widetilde{g}}}_{t})-V^\lambda_g(X^x_t,Y^{y,\xi^{\widetilde{g}}}_{t-})$ and $\Delta Y_t^{y,\xi^{\widetilde{g}}} =Y_{t}^{y,\xi^{\widetilde{g}}}-Y_{t-}^{y,\xi^{\widetilde{g}}}$. 
The second equality holds because 
    \begin{eqnarray*}
      &&   \mathbb{E}\Big[\int_0^{\tau_R\wedge T} e^{-\rho t} \Big(\mathcal{L}_x\,u(X^x_t,Y^{ y,\xi^g}_{t-})-\rho\, V^\lambda_g(X^x_t,Y^{ y,\xi^g}_{t-})\Big) \dd t \Big] \\
         &=& -\mathbb{E}\Big[\int_0^{\tau_R\wedge T} e^{-\rho t} \Big(\Big(\pi(X^x_t)-\rho \kappa\Big)Y_t^{ y,\xi^g} - \lambda  Y_t^{ y,\xi^g} \log\Big(Y_t^{ y,\xi^g}\Big) \Big)\dd t\Big],
    \end{eqnarray*} 
    as pure jumps could only possibly happen at time $0$,  $\xi^{\widetilde{g}}$ keeps $(X^x_t,Y_t^{y,\xi^{\widetilde{g}}})$ within $\overline{\mathcal{E}(g)}$ for $t >0$ and  \eqref{eq:HJB_g_1}. 
    The last inequality holds by the design of  \eqref{eq:updated_policy}. 
In particular,  given that $\widetilde{g}$ is increasing on $[0,\hat{x}_{\widetilde{g}}]$, we have $(x,y)\in \mathcal{S}( \widetilde{g})\cap \mathcal{E}(g)$ hold for all $x\in[g^{-1}(y),\widetilde{g}^{-1}(y)]$ and $y\in[0,1]$. Hence, for $(x,y)\in \mathcal{S}( \widetilde{g})\cap \mathcal{E}(g)$,
\begin{eqnarray*}
    \partial_y V^\lambda_g(x,y) =\partial_y V^\lambda_g(g^{-1}(y),y) + \int_{g^{-1}(y)}^x \partial_{xy} V^\lambda_g (u,y) {\rm d} u \leq 0,
\end{eqnarray*}
in which $\partial_{xy} V^\lambda_g (u,y) \leq 0$ for $(u,y)\in \mathcal{S}( \widetilde{g})\cap \mathcal{E}(g)$ due to the design of the algorithm.

In addition, it holds that 
\begin{eqnarray*}
    \Delta V^\lambda_g(X^x_t,Y^{y,\xi^{\widetilde{g}}}_{t}) &=&V^\lambda_g(X^x_t,Y^{y,\xi^{\widetilde{g}}}_{t})-V^\lambda_g(X^x_t,Y^{y,\xi^{\widetilde{g}}}_{t-}) \\
    &=& \int_0^{\Delta Y^{y,\xi^{\widetilde{g}}}_t} \partial_y V(X^x_t,Y^{y,\xi^{\widetilde{g}}}_{t-}+u) {\rm d} u \ge 0,
\end{eqnarray*}    
as $\Delta Y^{y,\xi^{\widetilde{g}}}_t \leq 0$.

When taking limits as $R \rightarrow\infty$ we have $\tau_R\wedge T\rightarrow T$, $\mathbb{P}$-a.s.  By standard properties of Brownian motion it is easy to prove that the integral terms in the last expression on the right-hand side of \eqref{eq:ito2} are uniformly bounded in $\mathbb{L}^2(\Omega, \mathbb{P})$, hence uniformly integrable. Moreover,
$V^\lambda_g$ (taking the form as in \eqref{eq:u_g}) has sublinear growth by straightforward calculation. Then we also take limits as $T\rightarrow \infty$ and it follows
that
\begin{eqnarray*}
    V^\lambda_{\widetilde{g}}(x,y) \ge V^\lambda_{g}(x,y).
\end{eqnarray*}


\subsection{Proof of Theorem \ref{thm:policy_convergence}}
\label{sec:AppE}

\emph{Step 1.} By the construction of the sequence in Algorithm \ref{alg1}, we have $ 0 \leq g_{k+1} \leq g_{k} \leq g_0$ for all $k\in \mathbb{N}_+$. Hence, by monotone convergence theorem, there exists a limit $\bar{g} \ge 0$ such that
\begin{eqnarray}
\bar{g}(x) := \lim_{k \rightarrow \infty} g_k(x). 
\end{eqnarray}
The main goal is to show that $\bar{g}(x)  = g_\lambda(x)$ for $x\in \mathbb{R}_+$.

\vskip10pt
\noindent \emph{Step 2.} Recall $V^\lambda_{g_k}(x,y)$ as the value function associated with strategy $\xi^{g_k}$. We use induction to prove the following iteratively:
\begin{itemize}
    \item[{\bf (a)}] $g_k$ is  strictly increasing on $[0,\hat{x}_{g_k}]$, 
    \item[{\bf (b)}] $g_k(0) {=} e^{-(1+
\frac{\kappa \rho}{\lambda})}$ for $k\in \mathbb{N}_+$,  
\item[{\bf (c)}] $- \alpha_- \big(\kappa  +\frac{\lambda}{\rho}  \log (g_k(x))+ \frac{\lambda}{\rho}\big) +\alpha_-H_\pi(x) - H_\pi^{\prime}(x)\cdot x \ge 0$ for $x\in [0,\hat{x}_{g_k}]$, namely, $g_k \ge g_\lambda$ on $[0,\hat{x}_{g_k}]$, 
{\item[{\bf (d)}] $\partial_{xy}^- u_{k}(x,g_k(x)) 
\leq 0$ for $x\in [0,\hat{x}_{g_k}]$.}
 \item[{\bf (e)}] $g_k\in {C}^1([0,\hat{x}_{g_k}])$, 
\end{itemize}
 The conditions {(a)-(c), (e)} hold for $k=0$ by the design of the initial reflection boundary $g_0$ (see Assumption \ref{ass:initial_policy}). Therefore, Theorem \ref{thm:property_g}  applies to $g_0$. Hence, as noted in the proof of Theorem \ref{thm:property_g},
 we have $u_0(x,y) = A_0(y) x^{\alpha_-} + H_\pi(x)\,y - \kappa y -\frac{\lambda}{\rho} y \log y$, with $A_0(y) = \int_{{g}_0(0)}^y  \frac{\kappa  +\frac{\lambda}{\rho}  \log (u)+ \frac{\lambda}{\rho}-H_\pi(g_0^{-1}(u))}{({g}_0^{-1}(u))^{\alpha_-}} \dd u$. Combined with condition {\bf (c)}, we have for $x\in [0,\hat{x}_{g_0}]$,
 \begin{eqnarray*}
      \partial_{yx}{u}_{0}(x,g_{0}(x)) =  \frac{\alpha_-}{x} \Big(\kappa  +\frac{\lambda}{\rho}  \log (g_{0}(x))+ \frac{\lambda}{\rho}-H_\pi(x) \Big)+ H_\pi'(x) \leq 0.
 \end{eqnarray*}
Therefore condition {\bf (d)} holds.
 
 Assume the above property holds for $k$. Then by Theorem \ref{thm:property_g}, 
 the following  HJB equation 
     \begin{eqnarray}
      &&(\mathcal{L}_x-\rho) u + \Big(\pi(x)-\rho \kappa \Big)y - \lambda y \log y                   =0\quad \text{on} \quad\mathcal{E}(g_k),\\
     && -u_y = 0, \quad \text{on}\quad \mathcal{S}(g_k).
  \end{eqnarray}
  has a unique ${C}^{1}(\mathbb{R}_+\times [0,1])\cap{C}^{2}\left(\overline{\mathcal{E}(g_k)}\right)$ solution, which we denoted as $u_k$. Also,  $V^\lambda_{g_k}(x,y) \equiv u_k(x,y)$.
 Then for $y \leq g_k(x)$ the value function $u_k$ satisfies the form:
 \begin{eqnarray}\label{eq:uk}
{u}_k(x,y) = A_k(y) x^{\alpha_-} + H_\pi(x) y - \kappa y -\frac{\lambda}{\rho} y \log y                  ,
\end{eqnarray}
with 
\begin{eqnarray}\label{eq:Ak}
 A_k(y) = \int_{{g}_k(0)}^y  \frac{\kappa  +\frac{\lambda}{\rho}  \log (u)+ \frac{\lambda}{\rho}-H_\pi(g_k^{-1}(u))}{({g}_k^{-1}(u))^{\alpha_-}} \dd u.
\end{eqnarray}
In this region, we also have,
\begin{eqnarray*}
    \partial_{y}{u}_{k}(x,y) &=& A_k'(y)x^{\alpha_-} + H_\pi(x) - k -\frac{\lambda}{\rho} -\frac{\lambda}{\rho}\log(y)\\
    &=& \frac{x^{\alpha_-}}{(g_k^{-1}(y)^{\alpha_-}} \Big(\kappa  +\frac{\lambda}{\rho}  \log y                  + \frac{\lambda}{\rho}-H_\pi(g_k^{-1}(y)) \Big)+ H_\pi(x) - k -\frac{\lambda}{\rho} -\frac{\lambda}{\rho}\log(y)\\
    &=& \frac{x^{\alpha_-}-(g_k^{-1}(y)^{\alpha_-}}{(g_k^{-1}(y)^{\alpha_-}} \Big(\kappa  +\frac{\lambda}{\rho}  \log y                  \\
    &&\qquad \qquad \qquad + \frac{\lambda}{\rho}-H_\pi(x) \Big) + \frac{x^{\alpha_-}}{(g_k^{-1}(y)^{\alpha_-}}\Big(H_\pi(x)- H_\pi(g_k^{-1}(y))\Big),
\end{eqnarray*}
and
\begin{eqnarray*}
    \partial_{yx}{u}_{k}(x,y) &=& \alpha_-\,A_k'(y)x^{\alpha_--1} + H_\pi'(x) \\
    &=& \frac{\alpha_-x^{\alpha_--1}}{(g_k^{-1}(y)^{\alpha_-}} \Big(\kappa  +\frac{\lambda}{\rho}  \log y                  + \frac{\lambda}{\rho}-H_\pi(g_k^{-1}(y)) \Big)+ H_\pi'(x).
\end{eqnarray*}


In the $k$th iteration there are two scenarios for each point $(x_0,g_k(x_0))$:  $\partial_{xy}^{-}u_k(x_0,g_k(x_0)) = 0$ or $\partial_{xy}^{-}u_k(x_0,g_k(x_0)) < 0$. When $\partial_{xy}^{-}u_k(x_0,g_k(x_0)) = 0$, we have $g_{k+1}(x_0) = g_{k}(x_0)$. Therefore, for both scenarios,
$\partial_{xy}^- u_k (x,y)=0$ leads to $y = {g}_{k+1}(x)$, or equivalently, we solve $y = {g}_{k+1}(x)\leq g_k(x)$ from $\partial_{xy}^- u_k \big(x,g_{k+1}(x)\big)=0$. Then we have
\begin{eqnarray}
0 &=& -\frac{\alpha_-x^{\alpha_-}}{(g_k^{-1}(g_{k+1}(x))^{\alpha_-}} \Big(\kappa  +\frac{\lambda}{\rho}  \log (g_{k+1}(x))+ \frac{\lambda}{\rho}-H_\pi(g_k^{-1}(g_{k+1}(x))) \Big)- H_\pi'(x)\, x\nonumber\\
&\leq& -\alpha_- \Big(\kappa  +\frac{\lambda}{\rho}  \log (g_{k+1}(x))+ \frac{\lambda}{\rho}-H_\pi(x) \Big)- H_\pi'(x)\, x.\label{eq:inner0}
\end{eqnarray}
The last inequality holds since for $y \ge g_\lambda(x)$,
\begin{eqnarray}\label{eq:inter1}
    \frac{k+\frac{\lambda}{\rho}\log(y)+\frac{\lambda}{\rho}-H_\pi(g_k^{-1}(g_{k+1}(x)))}{g_k^{-1}(g_{k+1}(x))^{\alpha_-}} \leq  \frac{k+\frac{\lambda}{\rho}\log(y)+\frac{\lambda}{\rho}-H_\pi(x)}{x^{\alpha_-}}.
\end{eqnarray}
Equation \eqref{eq:inter1} holds since $g_{k}^{-1} ({g}_{k+1}(x)) \leq x$ and 
\begin{eqnarray}
    &&\partial_x \Bigg(\frac{k+\frac{\lambda}{\rho}\log(y)+\frac{\lambda}{\rho}-H_\pi(x)}{x^{\alpha_-}}\Bigg) \nonumber\\
    &=& \frac{-H_\pi'(x)x^{\alpha_-}-\alpha_-\Big(k+\frac{\lambda}{\rho}\log(y)+\frac{\lambda}{\rho}-H_\pi(x)\Big)x^{\alpha_--1}}{x^{2\alpha_-} } \ge 0. \label{eq:inter2}
\end{eqnarray}
Finally, \eqref{eq:inter2} holds for $y \ge g_\lambda(x)$. 
Hence \eqref{eq:inner0} suggests that {\bf (c)} is satisfied in the $k+1$ iteration.

{
Next we show that $g_{k+1}(0) = e^{-(1+\frac{\kappa \rho}{\lambda})}$. 
For any points $x_0$ such that {\bf (c)} holds with equality sign in the  $k$th iteration, namely, 
$$
g_k(x_0) = \exp \left(\frac{\frac{-\alpha_-H_\pi(x_0)+H_\pi^{\prime}(x_0)x_0}{-\alpha_-}-(\kappa+\frac{\lambda}{\rho})\big)\rho}{\lambda}\right),
$$ 
it holds that $\partial_{xy}u_k^-(x_0,g_k(x_0)) = 0$ by simple calculation. 
Hence, $g_{k+1}(x_0)=g_{k}(x_0)$ and take $x_0=0$ to have $g_{k+1}(0) = g_k(0) = e^{-(1+\frac{\kappa \rho}{\lambda})}$.
}

We now show condition {\bf (a)} by contradiction. Suppose there exists $0<x_1<x_2$ such that $\bar{y} = g_{k+1}({x}_1) = g_{k+1}({x}_2)$. Then we have
\begin{eqnarray}
   0 &=&  x_1\partial_{xy}u_k(x_1,\bar{y}) -x_2\partial_{xy}u_k(x_2,\bar{y}) 
   \nonumber\\
   &=& 
    \frac{\alpha_-x_1^{\alpha_-}}{(g_k^{-1}(\bar y)^{\alpha_-}} \Big(\kappa  +\frac{\lambda}{\rho}  \log (\bar y)+ \frac{\lambda}{\rho}-H_\pi(g_k^{-1}(\bar y)) \Big)+ H_\pi'(x_1)x_1\nonumber\\
    && -\frac{\alpha_-x_2^{\alpha_-}}{(g_k^{-1}(\bar y)^{\alpha_-}} \Big(\kappa  +\frac{\lambda}{\rho}  \log (\bar y)+ \frac{\lambda}{\rho}-H_\pi(g_k^{-1}(\bar y)) \Big)- H_\pi'(x_2)x_2\nonumber\\
    &=& 
    \frac{\alpha_-(x_1^{\alpha_-}-x_x^{\alpha_-})}{(g_k^{-1}(\bar y)^{\alpha_-}} \Big(\kappa  +\frac{\lambda}{\rho}  \log (\bar y)+ \frac{\lambda}{\rho}-H_\pi(g_k^{-1}(\bar y)) \Big)+ H_\pi'(x_1)x_1-H_\pi'(x_2) x_2. \label{eq:proof-eqn-a}
\end{eqnarray}
Note that $x_1^{\alpha_-}>x_2^{\alpha_-}$ due to the facts that $x_1<x_2$ and $\alpha_-<0$. In addition, we have $x_1H_\pi'(x_1)<x_2H_\pi'(x_2)$ as $x\pi'(x)$ is strictly increasing according to Assumption \ref{ass:Pi}. Hence we get a contradiction to \eqref{eq:proof-eqn-a}. This suggests that $g_{k+1}$ is either strictly increasing or decreasing. Given that $0<g_\lambda \leq g_{k+1}$ on $[0,\hat{x}_{g_{k+1}}]$, then $g_{k+1}$ must be strictly increasing.

 {\bf (d)} holds since, according to {\bf (c)} 
 \begin{eqnarray*}
      \partial_{yx}{u}_{k+1}(x,g_{k+1}(x)) =  \frac{\alpha_-}{x} \Big(\kappa  +\frac{\lambda}{\rho}  \log (g_{k+1}(x))+ \frac{\lambda}{\rho}-H_\pi(x) \Big)+ H_\pi'(x) \leq 0.
 \end{eqnarray*}

Finally, given that $u_k\in {C}^{2}\left(\overline{\mathcal{E}(g_k)}\right)$, $\partial_{xy}u_k(x,g_{k+1}(x)) = 0$ for $x\in [0,\hat{x}_{g_{k+1}}]$ and $(x,g_{k+1}(x))\in \overline{\mathcal{E}(g_k)}$, we have $g_{k+1}\in {C}^1( [0,\hat{x}_{g_{k+1}}])$. Hence we have that  {\bf (e)} holds.

\vskip10pt
\noindent \emph{Step 3.} We study the property of the limiting function $\bar{g} \ge g_\lambda$ as $g_k \ge g_\lambda$ for all $k\in \mathbb{N}_+$. 

According to Step 2, {\bf (a)}-{\bf (e)} holds iteratively. Therefore it holds for all $k\in \mathbb{N}_+$ that
\begin{eqnarray}\label{eq:inter5}
\qquad    \frac{-\alpha_-x^{\alpha_-}}{(g_k^{-1}(g_{k+1}(x))^{\alpha_-}} \Big(\kappa  +\frac{\lambda}{\rho}  \log (g_{k+1}(x))+ \frac{\lambda}{\rho}-H_\pi(g_k^{-1}(g_{k+1}(x))) \Big)- H_\pi'(x)\, x=0.
\end{eqnarray}
For each $x\in\mathbb{R}_+$ fixed, we have $\lim_{k\rightarrow\infty}g_k^{-1}(g_{k+1}(x))=x$. 
Combined with \eqref{eq:inter5}, we find
\begin{eqnarray*}
    0 &=& -\alpha_- \Big(\kappa  +\frac{\lambda}{\rho}  \log (\bar g(x))+ \frac{\lambda}{\rho}-H_\pi(x) \Big)- H_\pi'(x)\, x.
\end{eqnarray*}
Therefore, $\bar g = g_\lambda$ for $x\in[0,x_{g_\lambda}]$. 
Hence we obtain \eqref{eq:policy_convergence1}.

In addition, according to Step 2, for $(x,y)\in \overline{\mathcal{E}(g_k)}$, we have $u_k$ follows \eqref{eq:uk} with $A_k$ defined in \eqref{eq:Ak}. As $g_k$ converges to $g_\lambda$ uniformly, we have $A_k$
converges to $A$, which is defined in \eqref{eq:A}. Hence $u_k$ converges to $u^\lambda$ in $\overline{\mathcal{E}(g_\lambda)}$. Hence $u_k = V^\lambda_{g_k}$ converges to $u^\lambda = V^\lambda$ on $\mathbb{R}_+\times[0,1]$.
Finally, \eqref{eq:policy_convergence2} holds.

\section{Conclusion}
\label{sec:conclusion}

We have developed a novel reinforcement learning framework for continuous-time real option problems through an entropy-regularized singular control formulation. By embedding the classical stopping problem into a probabilistic control structure, our approach enables exploration and learning in both model-based and model-free settings. We derived explicit analytical solutions for the regularized problem and showed that the optimal reflection boundary converges to the classical free-boundary in the vanishing entropy limit. Furthermore, we proposed and analyzed a policy iteration algorithm that learns the optimal boundary function from data, with theoretical guarantees on policy improvement and convergence. Our findings establish a rigorous connection between stochastic control theory and reinforcement learning, and suggest new avenues for data-driven approaches in high-dimensional or partially observed optimal stopping environments.

\vspace{5pt}

\noindent \textbf{Acknowledgments.} 
Supported  by the Deutsche Forschungsgemeinschaft (DFG, German Research Foundation) - Project-ID 317210226 - SFB 1283 and the NSF CAREER award DMS-2524465.

\bibliographystyle{siam}
\bibliography{main.bib,references}

\newpage

\end{document}

%% file: figures/g_update_v3.tex
\tikzset{every picture/.style={line width=0.75pt}} 

\begin{tikzpicture}[x=0.75pt,y=0.75pt,yscale=-1,xscale=1]

\draw    (122,161) -- (348.5,161) ;
\draw [shift={(350.5,161)}, rotate = 180] [color={rgb, 255:red, 0; green, 0; blue, 0 }  ][line width=0.75]    (10.93,-3.29) .. controls (6.95,-1.4) and (3.31,-0.3) .. (0,0) .. controls (3.31,0.3) and (6.95,1.4) .. (10.93,3.29)   ;
\draw    (122,161) -- (122.99,19) ;
\draw [shift={(123,17)}, rotate = 90.4] [color={rgb, 255:red, 0; green, 0; blue, 0 }  ][line width=0.75]    (10.93,-3.29) .. controls (6.95,-1.4) and (3.31,-0.3) .. (0,0) .. controls (3.31,0.3) and (6.95,1.4) .. (10.93,3.29)   ;
\draw    (247,50) -- (326,51) ;
\draw [color={rgb, 255:red, 74; green, 144; blue, 226 }  ,draw opacity=1 ][line width=0.75]    (197,50) -- (276,51) ;
\draw    (123,148) .. controls (197,143) and (233,108) .. (247,50) ;
\draw [color={rgb, 255:red, 208; green, 2; blue, 27 }  ,draw opacity=1 ]   (123,148) .. controls (197,146) and (262,109) .. (276,51) ;
\draw [color={rgb, 255:red, 208; green, 2; blue, 27 }  ,draw opacity=1 ]   (276,51) -- (326,51) ;
\draw [color={rgb, 255:red, 74; green, 144; blue, 226 }  ,draw opacity=1 ]   (191.5,102.5) -- (191.04,129) ;
\draw [shift={(191,131)}, rotate = 271.01] [color={rgb, 255:red, 74; green, 144; blue, 226 }  ,draw opacity=1 ][line width=0.75]    (10.93,-3.29) .. controls (6.95,-1.4) and (3.31,-0.3) .. (0,0) .. controls (3.31,0.3) and (6.95,1.4) .. (10.93,3.29)   ;
\draw [color={rgb, 255:red, 0; green, 0; blue, 0 }  ,draw opacity=1 ]   (238,76) -- (238,106) ;
\draw [shift={(238,108)}, rotate = 270] [color={rgb, 255:red, 0; green, 0; blue, 0 }  ,draw opacity=1 ][line width=0.75]    (10.93,-3.29) .. controls (6.95,-1.4) and (3.31,-0.3) .. (0,0) .. controls (3.31,0.3) and (6.95,1.4) .. (10.93,3.29)   ;
\draw [color={rgb, 255:red, 74; green, 144; blue, 226 }  ,draw opacity=1 ]   (123,148) .. controls (197,146) and (191,92) .. (197,50) ;

\draw (106,43) node [anchor=north west][inner sep=0.75pt]   [align=left] {1};
\draw (106,150) node [anchor=north west][inner sep=0.75pt]   [align=left] {0};
\draw (104,7.4) node [anchor=north west][inner sep=0.75pt]    {$y$};
\draw (357,152.4) node [anchor=north west][inner sep=0.75pt]    {$x$};
\draw (200,67.4) node [anchor=north west][inner sep=0.75pt]    {$g_{k}( x)$};
\draw (267,71.4) node [anchor=north west][inner sep=0.75pt]  [color={rgb, 255:red, 208; green, 2; blue, 27 }  ,opacity=1 ]  {$g_{k+1}( x)$};
\draw (153,90.4) node [anchor=north west][inner sep=0.75pt]  [color={rgb, 255:red, 74; green, 144; blue, 226 }  ,opacity=1 ]  {$g_{0}( x)$};

\end{tikzpicture}

%% file: main.bib
@article{dianetti.ferrari.xu.2024exploratory,
  title={Exploratory optimal stopping: {A} singular control formulation},
  author={Dianetti, Jodi and Ferrari, Giorgio and Xu, Renyuan},
  journal={arXiv preprint arXiv:2408.09335},
  year={2024}
}

@book{revuz.yor.2013continuous,
  title={Continuous {M}artingales and {B}rownian {M}otion},
  author={Revuz, Daniel and Yor, Marc},
  volume={293},
  year={2013},
  publisher={Springer Science \& Business Media}
}

@article{rao.chen.vemuri.wang.2004cumulative,
  title={Cumulative residual entropy: {A} new measure of information},
  author={Rao, Murali and Chen, Yunmei and Vemuri, Baba C and Wang, Fei},
  journal={IEEE Trans. Inform. Theory},
  volume={50},
  number={6},
  pages={1220--1228},
  year={2004},
  publisher={IEEE}
}

@book{pham2009continuous,
  title={Continuous-time {S}tochastic {C}ontrol and {O}ptimization with {F}inancial {A}pplications},
  author={Pham, Huy{\^e}n},
  volume={61},
  year={2009},
  publisher={Springer Science \& Business Media}
}

@book{fleming.soner2006,
  title={Controlled {M}arkov Processes and Viscosity Solutions},
  author={Fleming, Wendell H and Soner, Halil Mete}, 
  year={2006},
  publisher={Springer Science \& Business Media}
}

@book{DixitPindyck,
  title={Investment Under Uncertainty}, 
  author={Dixit, Robert K and Pindyck, Robert S},
  year={2012},
  publisher={Princeton University Press}
}

@article{angelis&ferrari&moriarty2019,
  title={A solvable two-dimensional degenerate singular stochastic control problem with nonconvex costs},
  author={De Angelis, Tiziano and Ferrari, Giorgio and Moriarty, John},
  journal={Math. Oper. Res.},
  volume={44},
  number={2},
  pages={512--531},
  year={2019},
  publisher={INFORMS}
}

@article{guo&tomecek2009,
  title={A class of singular control problems and the smooth fit principle},
  author={Guo, Xin and Tomecek, Pascal},
  journal={SIAM J. Control Optim.},
  volume={47},
  number={6},
  pages={3076--3099},
  year={2009},
  publisher={SIAM}
}

@article{K,
  title={Hedging and liquidation under transaction costs in currency markets},
  author={Kabanov, Yuri M},
  journal={Finance Stoch.},
  volume={3},
  number={2},
  pages={237--248},
  year={1999},
  publisher={Springer}
}

@article{T,
author = "Tarski, Alfred",
journal = "Pacific J. Math.",
number = "2",
pages = "285--309",
title = "A lattice-theoretical fixpoint theorem and its applications",
volume = "5",
year = "1955"
}


%% file: references.bib
@book{trigeorgis1996real,
  title={Real options: Managerial flexibility and strategy in resource allocation},
  author={Trigeorgis, Lenos},
  year={1996},
  publisher={MIT press}
}

@article{mezey2010real,
  title={Real options in resource economics},
  author={Mezey, Esther W and Conrad, Jon M},
  journal={Annu. Rev. Resour. Econ.},
  volume={2},
  number={1},
  pages={33--52},
  year={2010},
  publisher={Annual Reviews}
}

@article{xu2023decision,
  title={Decision making under costly sequential information acquisition: the paradigm of reversible and irreversible decisions},
  author={Xu, Renyuan and Zariphopoulou, Thaleia and Zhang, Luhao},
  journal={arXiv preprint arXiv:2401.00569},
  year={2023}
}

@article{alvarez2001adoption,
  title={Adoption of uncertain multi-stage technology projects: a real options approach},
  author={Alvarez, Luis HR and Stenbacka, Rune},
  journal={Journal of Mathematical Economics},
  volume={35},
  number={1},
  pages={71--97},
  year={2001},
  publisher={Elsevier}
}

@article{huisman2004strategic,
  title={Strategic technology adoption taking into account future technological improvements: A real options approach},
  author={Huisman, Kuno JM and Kort, Peter M},
  journal={European Journal of Operational Research},
  volume={159},
  number={3},
  pages={705--728},
  year={2004},
  publisher={Elsevier}
}

@article{bensoussan2019sequential,
  title={Sequential capacity expansion options},
  author={Bensoussan, Alain and Chevalier-Roignant, Beno{\^\i}t},
  journal={Operations research},
  volume={67},
  number={1},
  pages={33--57},
  year={2019},
  publisher={INFORMS}
}

@article{vintila2007real,
  title={Real options in capital budgeting. Pricing the option to delay and the option to abandon a project},
  author={Vintila, Nicoleta},
  journal={Theoretical and applied economics},
  volume={7},
  number={512},
  pages={47--58},
  year={2007},
  publisher={Asociatia Generala a Economistilor din Romania/Editura Economica}
}

@article{cruz2016assessing,
  title={Assessing the option to abandon an investment project by the binomial options pricing model.},
  author={Cruz Rambaud, Salvador and S{\'a}nchez P{\'e}rez, Ana Mar{\'\i}a},
  journal={Advances in decision sciences},
  year={2016}
}

@article{guo2005irreversible,
  title={Irreversible investment with regime shifts},
  author={Guo, Xin and Miao, Jianjun and Morellec, Erwan},
  journal={Journal of Economic Theory},
  volume={122},
  number={1},
  pages={37--59},
  year={2005},
  publisher={Elsevier}
}

@article{bulan2005real,
  title={Real options, irreversible investment and firm uncertainty: new evidence from US firms},
  author={Bulan, Laarni T},
  journal={Review of Financial Economics},
  volume={14},
  number={3-4},
  pages={255--279},
  year={2005},
  publisher={Elsevier}
}

@article{herrera2021optimal,
  title={Optimal stopping via randomized neural networks},
  author={Herrera, Calypso and Krach, Florian and Ruyssen, Pierre and Teichmann, Josef},
  journal={arXiv preprint arXiv:2104.13669},
  year={2021}
}

@book{stokey2008economics,
  title={The Economics of Inaction: Stochastic Control models with fixed costs},
  author={Stokey, Nancy L},
  year={2008},
  publisher={Princeton University Press}
}

@book{peskir.shiryav.2006optimal,
  title={Optimal stopping and free-boundary problems},
  author={Peskir, Goran and Shiryaev, Albert},
  year={2006},
  publisher={Springer}
}

@article{touzi.vieille.2002continuous,
  title={Continuous-time  {D}ynkin games with mixed strategies},
  author={Touzi, Nizar and Vieille, Nicolas},
  journal={SIAM J. Control Optim.},
  volume={41},
  number={4},
  pages={1073--1088},
  year={2002},
  publisher={SIAM}
}

@article{tran2024policy,
  title={Policy Iteration for exploratory {H}amilton--{J}acobi--{B}ellman equations},
  author={Tran, Hung Vinh and Wang, Zhenhua and Zhang, Yuming Paul},
  journal={arXiv preprint arXiv:2406.00612},
  year={2024}
}

@article{ma2024convergence,
  title={On Convergence and Rate of Convergence of Policy Improvement Algorithms},
  author={Ma, Jin and Wang, Gaozhan and Zhang, Jianfeng},
  journal={arXiv preprint arXiv:2406.10959},
  year={2024}
}

@inproceedings{thomas2016data,
  title={Data-efficient off-policy policy evaluation for reinforcement learning},
  author={Thomas, Philip and Brunskill, Emma},
  booktitle={International Conference on Machine Learning},
  pages={2139--2148},
  year={2016},
  organization={PMLR}
}

@article{duchi2015optimal,
  title={Optimal rates for zero-order convex optimization: The power of two function evaluations},
  author={Duchi, John C and Jordan, Michael I and Wainwright, Martin J and Wibisono, Andre},
  journal={IEEE Trans. Inf. Theory},
  volume={61},
  number={5},
  pages={2788--2806},
  year={2015},
  publisher={IEEE}
}

@article{christensen2023data,
  title={Data-driven rules for multidimensional reflection problems},
  author={Christensen, S{\"o}ren and Thomsen, Asbj{\o}rn Holk and Trottner, Lukas},
  journal={arXiv preprint arXiv:2311.06639},
  year={2023}
}

@article{bayer2021randomized,
  title={Randomized optimal stopping algorithms and their convergence analysis},
  author={Bayer, Christian and Belomestny, Denis and Hager, Paul and Pigato, Paolo and Schoenmakers, John},
  journal={SIAM Journal on Financial Mathematics},
  volume={12},
  number={3},
  pages={1201--1225},
  year={2021},
  publisher={SIAM}
}

@article{christensen2023nonparametric,
  title={Nonparametric learning for impulse control problems—Exploration vs. exploitation},
  author={Christensen, S{\"o}ren and Strauch, Claudia},
  journal={The Annals of Applied Probability},
  volume={33},
  number={2},
  pages={1569--1587},
  year={2023},
  publisher={Institute of Mathematical Statistics}
}

@article{christensen2024learning,
  title={Learning to reflect: A unifying approach for data-driven stochastic control strategies},
  author={Christensen, S{\"o}ren and Strauch, Claudia and Trottner, Lukas},
  journal={Bernoulli},
  volume={30},
  number={3},
  pages={2074--2101},
  year={2024},
  publisher={Bernoulli Society for Mathematical Statistics and Probability}
}

@article{christensen2023data2,
  title={Data-driven optimal stopping: A pure exploration analysis},
  author={Christensen, S{\"o}ren and Dexheimer, Niklas and Strauch, Claudia},
  journal={arXiv preprint arXiv:2312.05880},
  year={2023}
}

@book{sutton2018reinforcement,
  title={Reinforcement learning: An introduction},
  author={Sutton, Richard S and Barto, Andrew G},
  year={2018},
  publisher={MIT press}
}

@article{kumar2004numerical,
  title={A numerical method for solving singular stochastic control problems},
  author={Kumar, Sunil and Muthuraman, Kumar},
  journal={Oper. Res.},
  volume={52},
  number={4},
  pages={563--582},
  year={2004},
  publisher={INFORMS}
}

@article{ata2023singular,
  title={Singular Control of (Reflected) {B}rownian Motion: {A} Computational Method Suitable for Queueing Applications},
  author={Ata, Baris and Harrison, J Michael and Si, Nian},
  journal={arXiv preprint arXiv:2312.11823},
  year={2023}
}

@article{hambly2021policy,
  title={Policy gradient methods for the noisy linear quadratic regulator over a finite horizon},
  author={Hambly, Ben and Xu, Renyuan and Yang, Huining},
  journal={SIAM J. Control Optim.},
  volume={59},
  number={5},
  pages={3359--3391},
  year={2021},
  publisher={SIAM}
}

@article{agarwal2021theory,
  title={On the theory of policy gradient methods: Optimality, approximation, and distribution shift},
  author={Agarwal, Alekh and Kakade, Sham M and Lee, Jason D and Mahajan, Gaurav},
  journal={J. Mach. Learn. Res.},
  volume={22},
  number={98},
  pages={1--76},
  year={2021}
}

@inproceedings{fazel2018global,
  title={Global convergence of policy gradient methods for the linear quadratic regulator},
  author={Fazel, Maryam and Ge, Rong and Kakade, Sham and Mesbahi, Mehran},
  booktitle={International conference on machine learning},
  pages={1467--1476},
  year={2018},
  organization={PMLR}
}

@article{denkert2024control,
  title={Control randomisation approach for policy gradient
and application to reinforcement learning in optimal switching},
  author={Denkert, Robert   Pham, Huy{\^e}n and Warin, Xavier},
  journal={arXiv preprint arXiv:2404.17939},
  year={2024}
}

@article{soner2023stopping,
  title={Stopping times of boundaries: Relaxation and continuity},
  author={Soner, H Mete and Tissot-Daguette, Valentin},
  journal={arXiv preprint arXiv:2305.09766},
  year={2023}
}

@article{reppen2022neural,
  title={Neural optimal stopping boundary},
  author={Reppen, A Max and Soner, H Mete and Tissot-Daguette, Valentin},
  journal={arXiv preprint arXiv:2205.04595},
  year={2022}
}

@article{dai2024learning,
  title={Learning to Optimally Stop a Diffusion Process},
  author={Dai, Min and Sun, Yu and Xu, Zuo Quan and Zhou, Xun Yu},
  journal={arXiv preprint arXiv:2408.09242},
  year={2024}
}

@article{dong2024randomized,
  title={Randomized optimal stopping problem in continuous time and reinforcement learning algorithm},
  author={Dong, Yuchao},
  journal={SIAM J. Control Optim.},
  volume={62},
  number={3},
  pages={1590--1614},
  year={2024},
  publisher={SIAM}
}

@article{hare2019dealing,
  title={Dealing with sparse rewards in reinforcement learning},
  author={Hare, Joshua},
  journal={arXiv preprint arXiv:1910.09281},
  year={2019}
}

@article{devidze2022exploration,
  title={Exploration-guided reward shaping for reinforcement learning under sparse rewards},
  author={Devidze, Rati and Kamalaruban, Parameswaran and Singla, Adish},
  journal={Adv. Neural Inf. Process. Syst.},
  volume={35},
  pages={5829--5842},
  year={2022}
}

@article{menaldi1983some,
  title={On some cheap control problems for diffusion processes},
  author={Menaldi, Jos{\'e}-Luis and Robin, Maurice},
  journal={Trans. Amer. Math. Soc.},
  volume={278},
  number={2},
  pages={771--802},
  year={1983}
}

@article{bai2023reinforcement,
  title={Reinforcement Learning for optimal dividend problem under diffusion model},
  author={Bai, Lihua and Gamage, Thejani and Ma, Jin and Xie, Pengxu},
  journal={arXiv preprint arXiv:2309.10242},
  year={2023}
}

@article{huang2022convergence,
  title={Convergence of policy improvement for entropy-regularized stochastic control problems},
  author={Huang, Yu-Jui and Wang, Zhenhua and Zhou, Zhou},
  journal={arXiv preprint arXiv:2209.07059},
  year={2022}
}

@article{ferrari2015integral,
author = {Ferrari, Giorgio},
title = {{On an integral equation for the free-boundary of stochastic, irreversible investment problems}},
volume = {25},
journal = {Ann. Appl. Probab.},
number = {1},
publisher = {Institute of Mathematical Statistics},
pages = {150 -- 176},
year = {2015},
doi = {10.1214/13-AAP991},
URL = {https://doi.org/10.1214/13-AAP991}
}

@article{angelis2019solvable,
  title={A solvable two-dimensional degenerate singular stochastic control problem with nonconvex costs},
  author={De Angelis, Tiziano and Ferrari, Giorgio and Moriarty, John},
  journal={Math. Oper. Res.},
  volume={44},
  number={2},
  pages={512--531},
  year={2019},
  publisher={INFORMS}
}

@book{fleming2006controlled,
  title={Controlled Markov processes and viscosity solutions},
  author={Fleming, Wendell H and Soner, Halil Mete},
  volume={25},
  year={2006},
  publisher={Springer Science \& Business Media}
}

@article{ferrari2016irreversible,
  title={Irreversible investment under {L}{\'e}vy uncertainty: An equation for the optimal boundary},
  author={Ferrari, Giorgio and Salminen, Paavo},
  journal={Adv. in Appl. Probab.},
  volume={48},
  number={1},
  pages={298--314},
  year={2016},
  publisher={Cambridge University Press}
}
